\documentclass[12pt]{article}

\usepackage{amsmath, amssymb, amsfonts}
\usepackage{geometry}
\usepackage{fancyhdr}
\usepackage{hyperref}
\usepackage{graphicx, subcaption}
\usepackage{float}
\renewcommand{\mathbf}[1]{\boldsymbol{#1}}

\usepackage{helvet, caption}
\title{Learning Neural Pushforward Samplers for Distributions from Fokker-Planck Equations by Weak Adversarial Training}

\author{
  Andrew Qing He\\{\small Department of Mathematics}\\{\small Southern Methodist University}\\{\small \texttt{andrewho@smu.edu}} \and Wei Cai\\{\small Department of Mathematics}\\{\small Southern Methodist University}\\{\small \texttt{cai@smu.edu}}
}

\begin{document}

\maketitle

\begin{abstract}
This paper presents a new method for solving Fokker-Planck equations (FPE) by learning a neural sampler for the distribution given by the FPE via an adversarial training based on a weak formulation of the FPE where the adjoint operator of FPE acts on the test function. Such a weak formulation transforms the PDE solution problem into a Monte Carlo importance sampling problem where  the FPE solution-distribution is learned through a neural pushforward map, avoiding some of the limitations of direct PDE based methods. Moreover,  by using simple plane-wave test functions, derivatives on the test functions  can be explicitly computed. This approach produces a natural importance sampling strategy for the FPE solution distribution with probability conservation, from which the FPE solution can be easily constructed.
\end{abstract}


\vspace{1em}

\section{Introduction}
In recent years, significant attention has been directed towards solving differential and integral equations using deep learning techniques, such as Physics-Informed Neural Networks (PINNs) \cite{PINN}, the Deep Ritz method \cite{deepritz}, Weak Adversarial Networks (WAN) \cite{wan}, and others. A common element across these frameworks is the point-wise representation of the solution by a \textit{neural network}—with spatial or spatio-temporal coordinate as input and outputs as an estimate of the solution's value at the coordinates. This offers a natural and intuitive way to approximate  a PDE solution.

However, this approach is not always optimal, one of the main challenge is to find the optimal sampling strategy where the solution can be learned efficiently and accurately, especially in high dimensional cases. The nature of the solution may require a specific unknown representation or sampling. This paper consider a special class of PDEs, i.e.,  the solution as a probability density function (PDF) which is governed by a Fokker-Planck equation. Here, an approach by a pointwise evaluation of network is often suboptimal for two reasons:
\begin{enumerate}
    \item The solution may not be absolutely continuous with respect to the Lebesgue measure on $\mathbb{R}^n$ (e.g., a Dirac $\delta$ distribution), in which case a classical PDF does not exist.
    \item It is difficult to guarantee conservation of probability (a volume-preserving property), as exactly evaluating the integral of a neural network over an arbitrary domain is generally intractable.
\end{enumerate}

In such scenarios, a more natural strategy is to represent the solution distribution directly by a particle method as to be used in this work.  To obtain the desired distrbution, we propose representing the solution via a \textit{neural pushforward mapping} $F_{\boldsymbol{\vartheta}}$ (parameterized by ${\boldsymbol{\vartheta}}$), which transforms samples from a simple reference distribution (e.g., uniform or standard Gaussian) of potentially different dimensionality into the target solution distribution through a pushforward operation. This neural pushforward map will be learned by an adversial training of the weak form of the FPE obtained by integrating-by-parts twice such that the weak form can be viewed as statistical expectation with resepct to the FPE solution as the distribution.

A similar methodology has been employed in \cite{NPFP, WHF, XWanFP, TZhouFP}. However, these works typically employ a PINN-style loss that requires the existence and regularity of the PDF. This, in turn, constrains the network $F_{\boldsymbol{\vartheta}}$ to be invertible and Jacobian-friendly (e.g., Real NVP \cite{realNVP}, as used in normalizing flows), and incurs a high computational cost from calculating second-order derivatives of the density (which involves Jacobian determinant computations). In contrast, our approach is based on the \textit{weak} formulation of transport and Fokker-Planck equations. This allows us to employ more general deep neural networks for the mapping $F_{\boldsymbol{\vartheta}}$. We adopt the weak adversarial framework \cite{wan, weran}, which reformulates the equation in its weak form. A test function in this formulation is parameterized as an adversarial network (acting as a ``discriminator'', similar to GANs \cite{GAN}), and the pushforward mapping $F_{\boldsymbol{\vartheta}}$ is trained adversarially to satisfy the equation.

In our framework, $F_{\boldsymbol{\vartheta}}$ can be any neural network that maps from a convenient base distribution (of arbitrary dimension) to the target distribution space, with output dimension matching that of the target distribution. The trained network acts as an efficient sampler from the solution distribution. While the framework does not inherently provide an explicit expression for the probability density function, users who require one—and can ensure its existence—may still instantiate $F_{\boldsymbol{\vartheta}}$ as an invertible, Jacobian-friendly network like Real NVP. Nevertheless, this only serves as an optional choice within our framework, not a requirement.

The adversarial test function is selected in a special manner. Aligning with the philosophy of weak formulations, the test function should exhibit greater regularity and be more amenable to differential operators than the solution itself. In PDE theory, this implies the test functions should be smooth and exhibit suitable decay. But for numerically, this translates to a requirement for computational efficiency; the cost of differentiating the network should not substantially exceed the cost of evaluating it. To satisfy this condition, we select a specific, computationally efficient form for the test function: the plane wave function. This choice allows its partial derivatives to be computed explicitly and with the same complexity as evaluating the function itself, significantly enhancing the efficiency of the training process. It is worth noting that while plane waves have been employed as test functions in other works \cite{socmartnet, dfsocmartnet}, their treatment was fundamentally different as the requisite differential operators were not applied to them. Our method's explicit computation of these operators on the test function is a novel and innovative aspect of this framework.

\section{Methods}

\subsection{The Steady-State Fokker-Planck Equation}

We consider the steady-state Fokker-Planck equation (FPE) on $\mathbb{R}^{n}$, which describes the probability density function $\rho(\mathbf{x})$ of a stochastic process at equilibrium. The equation is given by:
\begin{equation} \label{eq:FPE}
\mathcal{L}^{\dagger} \rho(\mathbf{x}) = -\frac{1}{2}\sum _{i,j=1}^{n}\frac{\partial ^{2}}{\partial x_{i} \partial x_{j}}\left[ a_{ij}(\mathbf{x}) \rho(\mathbf{x})\right] +\sum _{i=1}^{n}\frac{\partial }{\partial x_{i}}\left[ b_{i}(\mathbf{x}) \rho(\mathbf{x})\right] = 0,
\end{equation}
where $\mathbf{x} = (x_1, \dots, x_n)^T \in \mathbb{R}^n$, $a_{ij}(\mathbf{x})$ are the components of the state-dependent diffusion tensor $\mathbf{A}(\mathbf{x})$, and $b_{i}(\mathbf{x})$ are the components of the drift vector $\mathbf{b}(\mathbf{x})$. We assume the diffusion tensor $\mathbf{A}(\mathbf{x})$ is such that the differential operator $\mathcal{L}^{\dagger}$ is elliptic (or possibly hypo-elliptic), ensuring the well-posedness of the problem. Furthermore, we seek a probability density function, which imposes the natural normalization constraint:
\begin{equation} \label{eq:normalization}
\int _{\mathbb{R}^{n}} \rho(\mathbf{x}) \, \mathrm{d} \mathbf{x} = 1, \quad \text{and} \quad \rho(\mathbf{x}) \geq 0.
\end{equation}

\subsubsection{Weak Formulation and Adversarial Formulation}

A key challenge in solving \eqref{eq:FPE} directly is that the solution $\rho$ may be a generalized function (e.g., a Dirac delta distribution) lacking a classical density. To circumvent this, we employ the weak formulation of the PDE. Multiplying \eqref{eq:FPE} by a smooth test function $f(\mathbf{x})$ and applying integration by parts twice (assuming appropriate decay conditions at infinity of $\rho$), we obtain the weak form:
\begin{equation} \label{eq:weak_form}
\int _{\mathbb{R}^{n}} \rho(\mathbf{x}) \left( -\frac{1}{2}\sum _{i,j=1}^{n} a_{ij}(\mathbf{x})\frac{\partial ^{2} f(\mathbf{x})}{\partial x_{i} \partial x_{j}} -\sum _{i=1}^{n} b_{i}(\mathbf{x})\frac{\partial f(\mathbf{x})}{\partial x_{i}}\right) \mathrm{d} \mathbf{x} = 0.
\end{equation}
This must hold for all test functions $f$ in a suitable class of smooth, bounded functions. Critically, this formulation does not require $\rho$ to be differentiable. The weak form can be interpreted as an expectation under the probability measure associated with $\rho$:
\begin{equation} \label{eq:expectation_form}
\mathbb{E}_{\mathbf{x} \sim \rho}\left[ \mathcal{L}f(\mathbf{x}) \right] = 0, \quad \text{where} \quad \mathcal{L}f(\mathbf{x}) = -\frac{1}{2}\sum _{i,j=1}^{n} a_{ij}(\mathbf{x})\frac{\partial ^{2} f}{\partial x_{i} \partial x_{j}} -\sum _{i=1}^{n} b_{i}(\mathbf{x})\frac{\partial f}{\partial x_{i}}.
\end{equation}
This expectation can be approximated efficiently via Monte Carlo sampling if one has access to samples from the distribution $\rho$.

\subsubsection{Neural Parametrization of the Solution and Test Functions}

Our method parametrizes two core components: the solution $\rho$ and the test functions $f$.

\paragraph{Parametrizing the Solution ($\rho$)}
Instead of representing the density $\rho$ directly, we represent the solution distribution implicitly via a pushforward map. We define a neural network $F_{{\boldsymbol{\vartheta}}}: \mathbb{R}^{d} \rightarrow \mathbb{R}^{n}$, parameterized by weights ${\boldsymbol{\vartheta}}$, that transforms samples from a simple base distribution (e.g., uniform or Gaussian) of dimension $d$ into samples from the target distribution $\rho$ in $\mathbb{R}^n$. The dimensions $d$ and $n$ need not be equal, providing flexibility in choosing the most convenient base distribution for sampling and optimization. If $\mathbf{r} \sim \mathcal{P}_{\text{base}}$, then $\mathbf{x} = F_{{\boldsymbol{\vartheta}}}(\mathbf{r})$ is a sample such that $\mathbf{x} \sim \rho$. This approach naturally enforces the normalization constraint and can represent distributions without a density.

\paragraph{Parametrizing the Test Functions ($f$)}
Following the rationale in the introduction, we choose test functions that are computationally efficient and facilitate easy computation of their derivatives. A highly effective choice is the \textit{plane-wave} function of the form:
\begin{equation*}
f(\mathbf{x}; \boldsymbol{\eta}^{(k)}) = \varphi\left( \mathbf{w}^{(k)} \cdot \mathbf{x} + b^{(k)} \right),
\end{equation*}
where $\boldsymbol{\eta}^{(k)} = \{\mathbf{w}^{(k)}, b^{(k)}\}$ are the parameters for the $k$-th test function, $\mathbf{w}^{(k)} \in \mathbb{R}^n$ is a weight vector, $b^{(k)} \in \mathbb{R}$ is a bias, and $\varphi: \mathbb{R} \rightarrow \mathbb{R}$ is a $C^2$ activation function. The derivatives of this function have a simple, closed-form expression. For a given test function $f$, the operator $\mathcal{L}f$ evaluates to:
\begin{align}
\mathcal{L}f(\mathbf{x}) = & -\left(\sum _{i=1}^{n} b_{i}(\mathbf{x}) w_{i}\right) \varphi '\left( \mathbf{w} \cdot \mathbf{x} + b \right) - \frac{1}{2}\left(\sum _{i,j=1}^{n} a_{ij}(\mathbf{x}) w_{i} w_{j}\right) \varphi ''\left( \mathbf{w} \cdot \mathbf{x} + b \right). \label{eq:L_f_general}
\end{align}
A common and effective choice is $\varphi(\cdot) = \sin(\cdot)$, leading to:
\begin{align}
\mathcal{L}f(\mathbf{x}) = & -\left(\sum _{i=1}^{n} b_{i}(\mathbf{x}) w_{i}\right)\cos\left( \mathbf{w} \cdot \mathbf{x} + b \right) + \frac{1}{2}\left(\sum _{i,j=1}^{n} a_{ij}(\mathbf{x}) w_{i} w_{j}\right)\sin\left( \mathbf{w} \cdot \mathbf{x} + b \right). \label{eq:L_f_sine}
\end{align}
This form is computationally efficient, particularly when the diffusion matrix $(a_{ij})_{n \times n}$ possesses favorable low-rank or sparsity structures. In such cases, the quadratic form $\sum_{i,j=1}^{n} a_{ij}(\mathbf{x}) w_i w_j$ can be evaluated in $\mathcal{O}(n)$ time, ensuring that the overall cost of computing $\mathcal{L}f$ remains comparable to that of evaluating $f$ itself. It is worth noting, however, that this efficiency is problem-dependent. For instance, when $(a_{ij})_{n \times n}$ is a constant multiple of the identity matrix—effectively reducing the second-order term to a Laplacian—the computation remains highly efficient. In the general case of a fully dense diffusion matrix, the double summation requires $\mathcal{O}(n^2)$ operations, leading to significantly increased computational cost in high dimensions. This scalability challenge is inherent to numerical methods for second-order partial differential equations and is not unique to the present approach.

\paragraph{Flexibility in Base Distribution Choice}
A key advantage of our approach is the flexibility in choosing the base distribution $\mathcal{P}_{\text{base}}$. Unlike methods that require the base and target distributions to have matching support or dimensionality, our framework allows $F_{{\boldsymbol{\vartheta}}}$ to map from any convenient $m$-dimensional base distribution to the $n$-dimensional target space. This flexibility enables:
\begin{itemize}
    \item Using higher-dimensional base distributions to provide richer representational capacity for complex target distributions
    \item Employing computationally efficient base distributions (e.g., uniform or standard normal) regardless of target complexity  
    \item Leveraging the universal approximation properties of neural networks to bridge arbitrary distributional gaps
\end{itemize}

\subsubsection{Adversarial Training Algorithm}

The training objective is to find parameters ${\boldsymbol{\vartheta}}$ for the pushforward network such that the weak formulation \eqref{eq:expectation_form} holds for a rich set of test functions. This is achieved by adversarially training the pushforward network $F_{{\boldsymbol{\vartheta}}}$ against a suite of $K$ parameterized test functions $\{f(\cdot; \boldsymbol{\eta}^{(k)})\}_{k=1}^K$. We denote $\boldsymbol{\eta}:= \{\boldsymbol{\eta}^{(k)}\}_{k=1}^K$ as the ensemble of all the parameters of the test functions.

The algorithm proceeds as follows:
\begin{enumerate}
    \item \textbf{Sampling:} Draw $M$ samples $\{\mathbf{r}^{(m)}\}_{m=1}^{M}$ from the base distribution $\mathcal{P}_{\text{base}}$ and push them through the network to generate samples $\{\mathbf{x}^{(m)} = F_{{\boldsymbol{\vartheta}}}(\mathbf{r}^{(m)})\}_{m=1}^{M}$.
    \item \textbf{Test Function Evaluation:} For each of the $K$ test functions with parameters $\boldsymbol{\eta}^{(k)} = \{\mathbf{w}^{(k)}, b^{(k)}\}$, compute the empirical mean of $\mathcal{L}f^{(k)}$ over the batch:
    \begin{equation*}
        \hat{\mathbb{E}}^{(k)}({\boldsymbol{\vartheta}}, \boldsymbol{\eta}^{(k)}) = \frac{1}{M} \sum_{m=1}^{M} \mathcal{L}f^{(k)}(\mathbf{x}^{(m)}).
    \end{equation*}
    \item \textbf{Loss Function:} The overall loss function measures the squared violation of the weak form across all $K$ test functions:
    \begin{equation} \label{eq:loss_function}
    {L}_{\text{total}}[ {\boldsymbol{\vartheta}} , \boldsymbol{\eta} ] = \frac{1}{K}\sum _{k=1}^{K} \left[ \hat{\mathbb{E}}^{(k)}({\boldsymbol{\vartheta}}, \boldsymbol{\eta}^{(k)}) \right]^2.
    \end{equation}
    \item \textbf{Adversarial Optimization:} We solve the following min-max optimization problem:
    \begin{equation*}
        \min_{{\boldsymbol{\vartheta}}} \max_{\boldsymbol{\eta}} {L}_{\text{total}}[ {\boldsymbol{\vartheta}} ,\boldsymbol{\eta} ].
    \end{equation*}
    In practice, this is done by alternating gradient updates:
    \begin{itemize}
        \item \textbf{Update ${\boldsymbol{\vartheta}}$ (Generator):} Take a gradient \textit{descent} step on ${\boldsymbol{\vartheta}}$ to minimize ${L}_{\text{total}}$, thus improving the pushforward map to satisfy the weak form for the current set of test functions.
        \item \textbf{Update $\boldsymbol{\eta}$ (Adversary):} Take a gradient \textit{ascent} step on each $\boldsymbol{\eta}$ to maximize ${L}_{\text{total}}$, thus refining the test functions to better detect violations of the weak form by the current pushforward map.
    \end{itemize}
\end{enumerate}
This adversarial training process ensures that the learned distribution $\rho$, represented by $F_{{\boldsymbol{\vartheta}}}$, converges to a solution that satisfies the Fokker-Planck equation in its weak form against a broad spectrum of test functions.

\subsection{Sampling and Integration with the Trained Network}

A key advantage of the pushforward representation is that the trained network $F_{\boldsymbol{\vartheta}}$ naturally provides an efficient sampler for the target distribution $\rho$. This enables straightforward Monte Carlo estimation of integrals with respect to the solution measure $\mu$ defined by $d\mu(x) = \rho(x)\,dx$. Specifically, for any integrable function $f: \mathbb{R}^n \to \mathbb{R}$, we can approximate:
\begin{equation}
\int_{\mathbb{R}^n} f(x)\,d\mu(x) \approx \frac{1}{M} \sum_{i=1}^M f(F_{\boldsymbol{\vartheta}}(r_i)),
\end{equation}
where $\{r_i\}_{i=1}^M$ are independent samples from the base distribution $\mathcal{P}_{\text{base}}$. This Monte Carlo estimator inherits the favorable convergence properties of importance sampling, with the pushforward map having been optimized to satisfy the Fokker-Planck equation in weak form.

\paragraph{Remark on Density Estimation}
While our framework naturally provides a sampler, it does not directly yield an explicit expression for the probability density function $\rho(x)$. In principle, if $F_{\boldsymbol{\vartheta}}$ is implemented as an invertible, Jacobian-friendly architecture (e.g., Real NVP~\cite{realNVP}), one could recover the density via the change-of-variables formula:
\begin{equation}
\rho(x) = \tilde{\rho}(F_{\boldsymbol{\vartheta}}^{-1}(x)) \left|\det J_{F_{\boldsymbol{\vartheta}}}(x)\right|^{-1},
\end{equation}
where $\tilde{\rho}$ is the density of the base distribution and $J_{F_{\boldsymbol{\vartheta}}}$ is the Jacobian of $F_{\boldsymbol{\vartheta}}$. However, we do not recommend this approach for two reasons:

\begin{enumerate}
    \item As a weak formulation method focused on satisfying integrated conditions rather than pointwise density values, our framework typically does not learn density profiles with high pointwise accuracy (as evidenced by the ring distribution example in Section~\ref{sec:ring}, where the radial density width deviates from the analytical prediction despite excellent agreement in the support and moment statistics).
    
    \item For complex target distributions, using a higher-dimensional base distribution ($d > n$) often provides superior representational capacity and training stability. This dimensional mismatch precludes the use of invertible architectures, which require $d = n$.
\end{enumerate}

Therefore, we emphasize that the primary output of our method is a sampler for computing expectations and generating realizations from the solution distribution, rather than an explicit density function. Users requiring explicit density estimates should either: (i) employ kernel density estimation on the generated samples, or (ii) restrict to invertible architectures with $d = n$, accepting the associated computational overhead and representational constraints.

\subsection{Extension to Time-Dependent Fokker-Planck Equations}

We now extend the framework to handle time-dependent Fokker-Planck equations with initial conditions. Consider the initial value problem on $[0, T] \times \mathbb{R}^n$:
\begin{equation} \label{eq:time_FPE}
\frac{\partial \rho}{\partial t} - \frac{1}{2}\sum_{i,j=1}^{n}\frac{\partial^2}{\partial x_i \partial x_j}[a_{ij}(t,\mathbf{x}) \rho] + \sum_{i=1}^{n}\frac{\partial}{\partial x_i}[b_i(t,\mathbf{x}) \rho] = 0,
\end{equation}
with initial condition $\rho(0, \mathbf{x}) = \rho_0(\mathbf{x})$, where $\int_{\mathbb{R}^n} \rho(t, \mathbf{x}) \, d\mathbf{x} = 1$ for all $t \in [0,T]$.

\subsubsection{Weak Formulation}

Multiplying by a test function $f(t, \mathbf{x})$ and integrating by parts in both space and time yields the weak formulation:
\begin{equation} \label{eq:time_weak_form}
\mathbb{E}_{\rho(T,\cdot)}[f(T,\cdot)] - \mathbb{E}_{\rho_0}[f(0,\cdot)] - \int_0^T \mathbb{E}_{\mathbf{x} \sim \rho(t,\cdot)}\left[\frac{\partial f}{\partial t} + \mathcal{L}f(t, \mathbf{x})\right] dt = 0,
\end{equation}
where $\mathcal{L}$ is the spatial Fokker-Planck operator applied to the test function:
\begin{equation}
\mathcal{L}f = -\frac{1}{2}\sum_{i,j=1}^{n} a_{ij}\frac{\partial^2 f}{\partial x_i \partial x_j} - \sum_{j=1}^{n} b_j\frac{\partial f}{\partial x_j}.
\end{equation}

This formulation decomposes the problem into three terms: a terminal expectation at time $T$, an initial condition expectation at time $0$, and an integral of expectations over the time interval $[0, T]$.

\subsubsection{Time-Parameterized Pushforward Map}

For time-dependent problems, we extend the pushforward map to depend on both time and the initial distribution. The key insight is to construct a network that naturally satisfies the initial condition while learning the temporal evolution. We parametrize:
\begin{equation}
F_{\boldsymbol{\vartheta}}(t, x_0, \mathbf{r}) = x_0 + \sqrt{t} \, \tilde{F}_{\boldsymbol{\vartheta}}(t, \mathbf{r}),
\end{equation}
where:
\begin{itemize}
\item $x_0 \sim \rho_0$ is sampled from the initial distribution $\rho_0$
\item $\mathbf{r} \sim \mathcal{P}_{\text{base}}$ is sampled from a simple base distribution (e.g., standard Gaussian)
\item $\tilde{F}_{\boldsymbol{\vartheta}}: \mathbb{R}^{1+d} \to \mathbb{R}^n$ is a neural network mapping $(t, \mathbf{r})$ to displacement vectors
\end{itemize}

This construction has several important properties:
\begin{enumerate}
\item \textbf{Initial condition enforcement}: At $t=0$, we have $F_{\boldsymbol{\vartheta}}(0, x_0, \mathbf{r}) = x_0$, automatically satisfying $\rho(0, \cdot) = \rho_0$
\item \textbf{Brownian scaling}: The $\sqrt{t}$ factor matches the characteristic scaling of diffusion processes, where mean squared displacement grows as $\mathcal{O}(t)$, so typical displacement scales as $\mathcal{O}(\sqrt{t})$
\item \textbf{Smooth initialization}: The network $\tilde{F}_{\boldsymbol{\vartheta}}$ does not need to learn the singularity at $t=0$ that would arise from a Brownian motion starting from rest
\end{enumerate}

During training, for each sample, we draw $t \sim \mathcal{U}(\epsilon, T)$ (with small $\epsilon > 0$ to avoid numerical issues), $x_0 \sim \rho_0$, and $\mathbf{r} \sim \mathcal{P}_{\text{base}}$, then compute $\mathbf{x}(t) = F_{\boldsymbol{\vartheta}}(t, x_0, \mathbf{r})$.

\subsubsection{Test Functions with Temporal Component}

We extend the plane-wave test functions to include temporal dependence:
\begin{equation}
f^{(k)}(t, \mathbf{x}) = \sin\left(\sum_{i=1}^{n} w_i^{(k)} x_i + \kappa^{(k)} t + b^{(k)}\right),
\end{equation}
where $\{\mathbf{w}^{(k)} \in \mathbb{R}^n, \kappa^{(k)} \in \mathbb{R}, b^{(k)} \in \mathbb{R}\}_{k=1}^K$ are trainable parameters. The temporal frequency $\kappa^{(k)}$ allows the test function to capture time-dependent features of the solution.

The derivatives required for the weak form have explicit expressions:
\begin{align}
\frac{\partial f^{(k)}}{\partial t} &= \kappa^{(k)} \cos\left(\mathbf{w}^{(k)} \cdot \mathbf{x} + \kappa^{(k)} t + b^{(k)}\right), \\
\frac{\partial f^{(k)}}{\partial x_i} &= w_i^{(k)} \cos\left(\mathbf{w}^{(k)} \cdot \mathbf{x} + \kappa^{(k)} t + b^{(k)}\right), \\
\frac{\partial^2 f^{(k)}}{\partial x_i \partial x_j} &= -w_i^{(k)} w_j^{(k)} \sin\left(\mathbf{w}^{(k)} \cdot \mathbf{x} + \kappa^{(k)} t + b^{(k)}\right).
\end{align}

These explicit forms enable efficient computation of $\mathcal{L}f^{(k)}$ without automatic differentiation.

\subsubsection{Monte Carlo Discretization and Loss Function}

The weak formulation \eqref{eq:time_weak_form} is discretized using Monte Carlo sampling. We approximate the three terms separately:

\paragraph{Terminal term ($t=T$):} Sample $M_T$ points with $t = T$:
\begin{equation}
\hat{E}_T^{(k)} = \frac{1}{M_T} \sum_{m=1}^{M_T} f^{(k)}\left(T, F_{\boldsymbol{\vartheta}}(T, x_{0,T}^{(m)}, r_T^{(m)})\right),
\end{equation}
where $x_{0,T}^{(m)} \sim \rho_0$ and $r_T^{(m)} \sim \mathcal{P}_{\text{base}}$.

\paragraph{Initial condition term ($t=0$):} Sample $M_0$ points directly from $\rho_0$:
\begin{equation}
\hat{E}_0^{(k)} = \frac{1}{M_0} \sum_{m=1}^{M_0} f^{(k)}(0, x_{0,0}^{(m)}),
\end{equation}
where $x_{0,0}^{(m)} \sim \rho_0$.

\paragraph{Interior integral term:} Sample $M$ space-time points uniformly over $[\epsilon, T]$:
\begin{equation}
\hat{E}^{(k)} = \frac{T - \epsilon}{M} \sum_{m=1}^{M} \left[\frac{\partial f^{(k)}}{\partial t} + \mathcal{L}f^{(k)}\right]\left(t^{(m)}, F_{\boldsymbol{\vartheta}}(t^{(m)}, x_0^{(m)}, r^{(m)})\right),
\end{equation}
where $t^{(m)} \sim \mathcal{U}(\epsilon, T)$, $x_0^{(m)} \sim \rho_0$, and $r^{(m)} \sim \mathcal{P}_{\text{base}}$. The factor $(T - \epsilon)$ accounts for the measure of the time domain.

\paragraph{Combined loss:} For each test function $k$, the residual of the weak form is:
\begin{equation}
R^{(k)} = \hat{E}_T^{(k)} - \hat{E}_0^{(k)} - \hat{E}^{(k)}.
\end{equation}
The total loss aggregates over all $K$ test functions:
\begin{equation}
\mathcal{L}_{\text{total}}[\boldsymbol{\vartheta}, \{\boldsymbol{\eta}^{(k)}\}] = \frac{1}{K} \sum_{k=1}^{K} \left(R^{(k)}\right)^2,
\end{equation}
where $\boldsymbol{\eta}^{(k)} = \{\mathbf{w}^{(k)}, \kappa^{(k)}, b^{(k)}\}$ are the test function parameters.

\paragraph{Adversarial optimization:} We solve:
\begin{equation}
\min_{\boldsymbol{\vartheta}} \max_{\{\boldsymbol{\eta}^{(k)}\}} \mathcal{L}_{\text{total}}[\boldsymbol{\vartheta}, \{\boldsymbol{\eta}^{(k)}\}],
\end{equation}
by alternating gradient descent on $\boldsymbol{\vartheta}$ (minimizing) and gradient ascent on $\{\boldsymbol{\eta}^{(k)}\}$ (maximizing). This adversarial training ensures the learned distribution satisfies the weak formulation against a broad and adaptive set of test functions.

\section{Results}

\subsection{A One-Dimensional Steady-State Problem}

To validate the proposed methodology, we first apply it to a canonical one-dimensional problem with a known analytical solution. This allows for a precise quantitative evaluation of the algorithm's accuracy and convergence properties.

\subsubsection{Problem Specification}
We consider the steady-state Fokker-Planck equation \eqref{eq:FPE} on $\mathbb{R}$ with the following parameters:
\begin{itemize}
    \item Dimension: $n=1$
    \item Drift coefficient: $b(x) = -\theta (x - \mu)$
    \item Diffusion coefficient: $a(x) = \sigma^2/2$ ($\sigma$ is a constant)
\end{itemize}
The corresponding PDE is therefore:
\begin{equation} \label{eq:1d_fpe}
-\frac{\sigma}{2} \frac{\partial ^{2} \rho}{\partial x^{2}} + \frac{\partial }{\partial x}\left[ (-\theta (x-\mu)) \rho\right] = 0.
\end{equation}
For this test, we set the parameters to $\sigma = \sqrt{2}$, $\theta = 1$, and $\mu = 2$, simplifying the equation to:
\begin{equation}
-\frac{\partial ^{2} \rho}{\partial x^{2}} + \frac{\partial (( x-2) \rho)}{\partial x} = 0.
\end{equation}
The analytical solution to this equation, satisfying the normalization constraint of a probability density function, is the Gaussian distribution $\mathcal{N}(\mu, \sigma^2/2\theta)$:
\begin{equation}
\rho_{\text{true}}( x) = \frac{1}{\sqrt{2\pi}}\exp\left({-\frac{( x-2)^{2}}{2}}\right).
\end{equation}
Thus, the target solution is a normal distribution with mean $\mu_{\text{true}} = 2$ and standard deviation $\sigma_{\text{true}} = 1$.

\subsubsection{Training Configuration}
The neural pushforward mapping $F_{{\boldsymbol{\vartheta}}}$ was implemented as a deep neural network (DNN) with architecture $(1, 10, 10, 1)$, using smooth activation functions (e.g., SiLU or Tanh) to ensure differentiability. This network maps samples from a simple base distribution to samples approximating the target distribution. In this experiment, the base distribution was chosen as the uniform distribution on $[0, 1]$, i.e., $r^{(m)} \sim \mathcal{U}(0,1)$.

The adversarial test functions were parameterized as $f^{(k)}(x) = \sin(w^{(k)} x + b^{(k)})$. A bank of $K=100$ such test functions was constructed. Notably, for this initial proof-of-concept, the parameters $\eta = \{w^{(k)}, b^{(k)}\}$ were not trained adversarially. Instead, the weights $w^{(k)}$ were initialized by sampling from a standard normal distribution $\mathcal{N}(0, 1)$, and the biases $b^{(k)}$ were sampled from $\mathcal{U}(0, 2\pi)$; these parameters remained fixed throughout the training process. This simplification tests the method's ability to find a solution that satisfies the weak form against a fixed, random set of modes. The batch size was set to $M=1000$ samples.

The model was trained by minimizing the loss function $\mathcal{L}_{\text{total}}$ (defined as \eqref{eq:loss_function}) using the Adam optimizer, performing gradient descent \textit{only} on the parameters ${\boldsymbol{\vartheta}}$ of the pushforward network $F_{{\boldsymbol{\vartheta}}}$.

\subsubsection{Results and Performance}
\begin{figure}[htbp]
    \centering
    \begin{subfigure}[t]{0.45\textwidth}
        \centering
        \includegraphics[width=\textwidth]{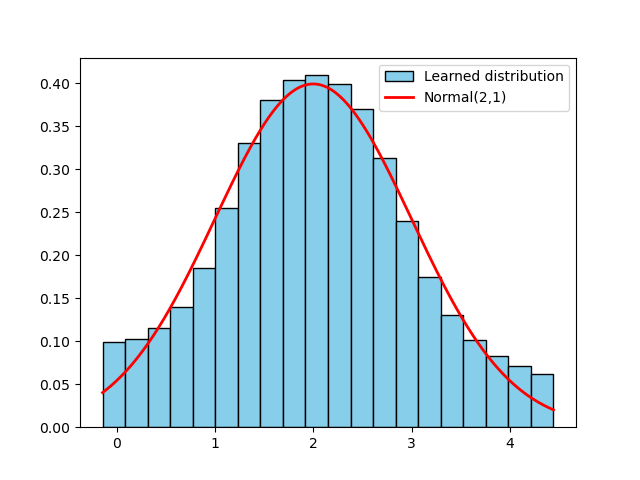}
        \caption{Comparison of the learned and true probability density functions.}
        \label{fig:1d_result:a}
    \end{subfigure}
    \hspace{1.5em}
    \begin{subfigure}[t]{0.45\textwidth}
        \centering
        \includegraphics[width=\textwidth]{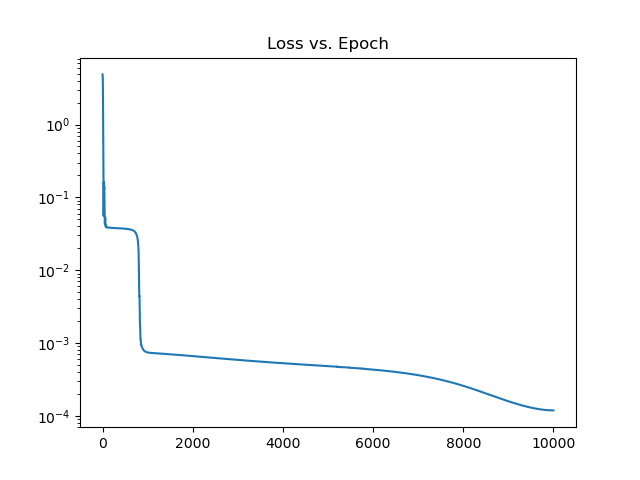}
        \caption{The convergence of the loss function over training steps (epochs).}
        \label{fig:1d_result:b}
    \end{subfigure}
    \caption{Numerical results for the one-dimensional Fokker-Planck equation. (a) The pushforward network $F_{{\boldsymbol{\vartheta}}}$ successfully learns to transform a uniform distribution into the target Gaussian distribution $\mathcal{N}(2,1)$, as the learned PDF shows excellent agreement with the analytical solution. (b) The training dynamics show a rapid and stable decrease in the loss value $\mathcal{L}_{\text{total}}$, indicating the model is effectively minimizing the violation of the PDE's weak form.}
    \label{fig:1d_result}
\end{figure}

The training converged successfully, with the loss decreasing from an initial value of $4.875$ to a final value of $0.000119$. Figure~\ref{fig:1d_result}(a) shows near-perfect alignment between the learned and analytical PDFs. The empirical mean and standard deviation of the learned distribution are $\mu_{\text{learned}} = 2.018553$ and $\sigma_{\text{learned}} = 0.974400$, showing excellent agreement with the true values ($\mu_{\text{true}}=2$ and $\sigma_{\text{true}}=1$) with relative errors below $1\%$ and $2.6\%$, respectively.

\subsection{A One-Dimensional Double-Peak Distribution}

We apply the method to a bimodal distribution arising from a double-well potential, showcasing the framework's ability to handle complex, multimodal probability distributions.

\subsubsection{Problem Specification}
We consider a steady-state Fokker-Planck equation with a double-well potential $V(x) = (x^2 - 1)^2$, which creates two symmetric minima at $x = \pm 1$. The corresponding drift and diffusion coefficients are:
\begin{itemize}
    \item Dimension: $n=1$
    \item Drift coefficient: $b(x) = -\frac{dV}{dx} = -4x^3 + 4x = -4x(x^2 - 1)$
    \item Diffusion coefficient: $a(x) = \sigma^2$ (constant)
\end{itemize}

For this experiment, we set $\sigma = 0.6$. The resulting Fokker-Planck equation is:
\begin{equation}
-\frac{\sigma^2}{2} \frac{\partial^2 \rho}{\partial x^2} + \frac{\partial}{\partial x}[(-4x^3 + 4x) \rho] = 0.
\end{equation}

The steady-state solution corresponds to the Boltzmann distribution:
\begin{equation}
\rho_{\text{true}}(x) \propto \exp\left(-\frac{2V(x)}{\sigma^2}\right) = \exp\left(-\frac{2(x^2-1)^2}{\sigma^2}\right),
\end{equation}
which exhibits two pronounced peaks near $x = \pm 1$ when the noise level $\sigma$ is sufficiently small.

\subsubsection{Training Configuration}
The neural pushforward mapping was implemented as a deeper network with architecture $(8, 100, 1)$ using Tanh activation functions, taking 8-dimensional input from a multivariate standard normal base distribution: 
$
d = 8,\quad \mathcal{P}_{\text{base}} = \mathcal{N}(0, I_{d\times d}).
$
This higher-dimensional input space provides the network with greater flexibility to model the complex bimodal structure.

The adversarial test functions maintained the same plane-wave parameterization $f^{(k)}(x) = \sin(w^{(k)} x + b^{(k)})$, with $K = 100$ test functions. Unlike the previous example, we implemented full adversarial training, alternating between:
\begin{itemize}
    \item Generator updates: minimizing the loss with respect to the pushforward network parameters
    \item Adversary updates: maximizing the loss with respect to test function parameters, performed 5 times per generator update every 10 epochs
\end{itemize}

The model was trained for 16,000 epochs with a batch size of 1,000 samples, using Adam optimizer with learning rate $2 \times 10^{-3}$ for the generator and SGD with learning rate $1 \times 10^{-2}$ for the test functions.

\subsubsection{Results and Performance}
\begin{figure}[htbp]
    \centering
    \begin{subfigure}[b]{0.45\textwidth}
        \centering
        \includegraphics[width=\textwidth]{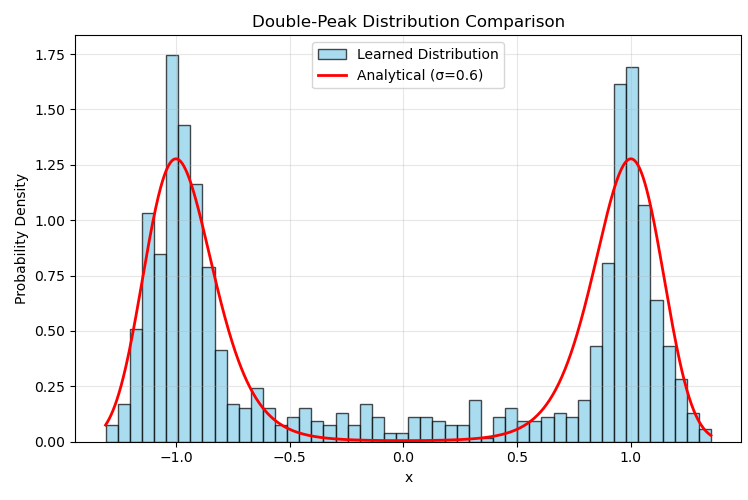}
        \caption{Comparison of the learned and analytical double-peak distributions.}
        \label{fig:double_peak:a}
    \end{subfigure}
    \hfill
    \begin{subfigure}[b]{0.45\textwidth}
        \centering
        \includegraphics[width=\textwidth]{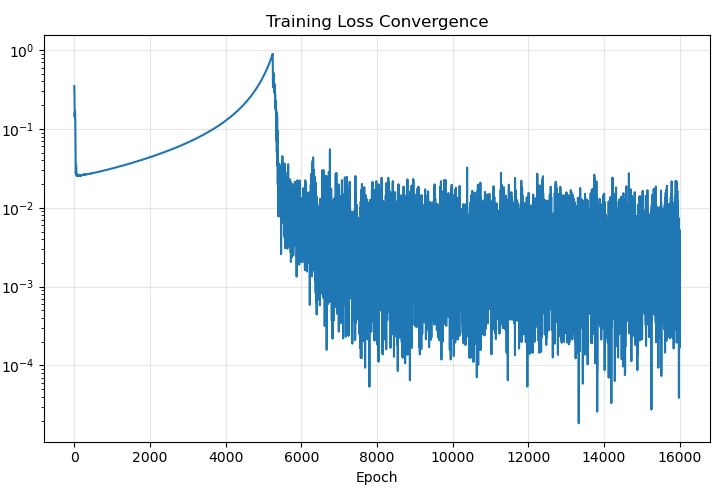}
        \caption{Training loss convergence for the double-peak problem.}
        \label{fig:double_peak:b}
    \end{subfigure}
    \caption{Numerical results for the one-dimensional double-peak Fokker-Planck equation. (a) The learned distribution successfully captures the bimodal structure with peaks at $x \approx \pm 1$, showing excellent agreement with the analytical Boltzmann distribution. (b) The adversarial training exhibits stable convergence with the loss decreasing over several orders of magnitude.}
    \label{fig:double_peak}
\end{figure}

The training demonstrated stable convergence with the loss function decreasing steadily over multiple orders of magnitude, as shown in Figure~\ref{fig:double_peak}(b). The neural pushforward mapping successfully transformed samples from a high-dimensional Gaussian base distribution into the target bimodal distribution, with clear double peaks at approximately $x = \pm 1$ matching the analytical Boltzmann distribution (Figure~\ref{fig:double_peak}(a)).

\subsection{A Two-Dimensional Double-Peak Problem}

We extend the framework to a two-dimensional problem with two distinct probability peaks to demonstrate the method's capability in handling multimodal distributions in higher dimensions.

\subsubsection{Problem Specification}
We consider the steady-state Fokker-Planck equation \eqref{eq:FPE} on $\mathbb{R}^2$ with the following parameters:
\begin{itemize}
    \item Dimension: $n=2$
    \item Drift coefficient: $\mathbf{b}(\mathbf{x}) = -\nabla V(\mathbf{x})$, where $V(x_1, x_2) = [(x_1-1)^2 + (x_2-1)^2][(x_1+1)^2 + (x_2+1)^2]$
    \item Diffusion coefficient: $a_{ij}(\mathbf{x}) = \sigma^2 \delta_{ij}/2$ (isotropic constant diffusion)
\end{itemize}

The potential function $V(\mathbf{x})$ creates a double-well landscape with minima located at $(-1, -1)$ and $(1, 1)$. The drift is computed as:
\begin{align}
b_1(\mathbf{x}) &= -2(x_1-1)[(x_1+1)^2 + (x_2+1)^2] - 2(x_1+1)[(x_1-1)^2 + (x_2-1)^2], \\
b_2(\mathbf{x}) &= -2(x_2-1)[(x_1+1)^2 + (x_2+1)^2] - 2(x_2+1)[(x_1-1)^2 + (x_2-1)^2].
\end{align}

For this test, we set $\sigma = 0.4$. The analytical steady-state solution follows the Boltzmann distribution:
\begin{equation} \label{ring-rho-true}
\rho_{\text{true}}(\mathbf{x}) \propto \exp\left(-\frac{2V(\mathbf{x})}{\sigma^2}\right),
\end{equation}
which exhibits two prominent peaks centered at the potential minima.

\subsubsection{Training Configuration}
The neural pushforward mapping $F_{\boldsymbol{\vartheta}}$ was implemented as a compact fully connected network with layer sizes $(16, 32, 32, 2)$, using Tanh activation functions. The input dimension of 16 provides a higher-dimensional base space to facilitate better mixing and exploration during training. The base distribution was chosen as a 16-dimensional standard Gaussian, $\mathbf{r} \sim \mathcal{N}(\mathbf{0}, \mathbf{I}_{16})$.

The adversarial test functions were parameterized as $f^{(k)}(\mathbf{x}) = \sin(\mathbf{w}^{(k)} \cdot \mathbf{x} + b^{(k)})$, with $K=200$ test functions to adequately capture the two-dimensional solution space. Both the test function parameters and the pushforward network were trained using their respective optimizers: stochastic gradient descent with learning rate $10^{-2}$ for the test functions, and Adam optimizer with learning rate $10^{-2}$ for the generator network. The batch size was set to $M=1000$ samples.

The model was trained for 100,000 epochs. During training, the test functions were updated adversarially, with 5 gradient ascent steps per generator update every 10 epochs, and a single update otherwise. This asymmetric update schedule helps maintain balance in the adversarial dynamics while ensuring the test functions remain sufficiently expressive to detect violations of the weak form.

\subsubsection{Results and Performance}

\begin{figure}[htbp]
    \centering
    \includegraphics[width=\textwidth]{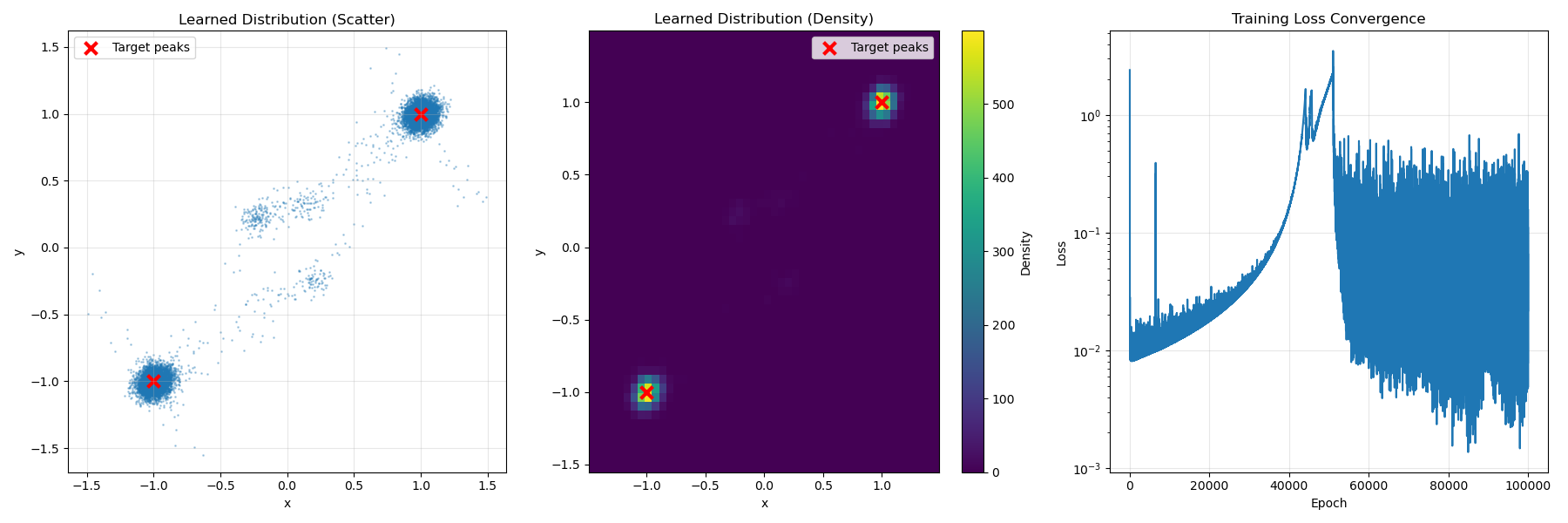}
    \caption{Training results for the two-dimensional double-peak Fokker-Planck equation. Left: Scatter plot of samples generated by the learned pushforward map, showing clear clustering at the target peaks $(-1,-1)$ and $(1,1)$ (marked with red crosses). Middle: Density heatmap revealing two distinct high-probability regions at the correct locations. Right: Training loss convergence over 100,000 epochs, demonstrating initial rapid decrease followed by a stable regime, with a transient instability around epoch 50,000--60,000.}
    \label{fig:2d_result_training}
\end{figure}

\begin{figure}[htbp]
    \centering
    \includegraphics[width=\textwidth]{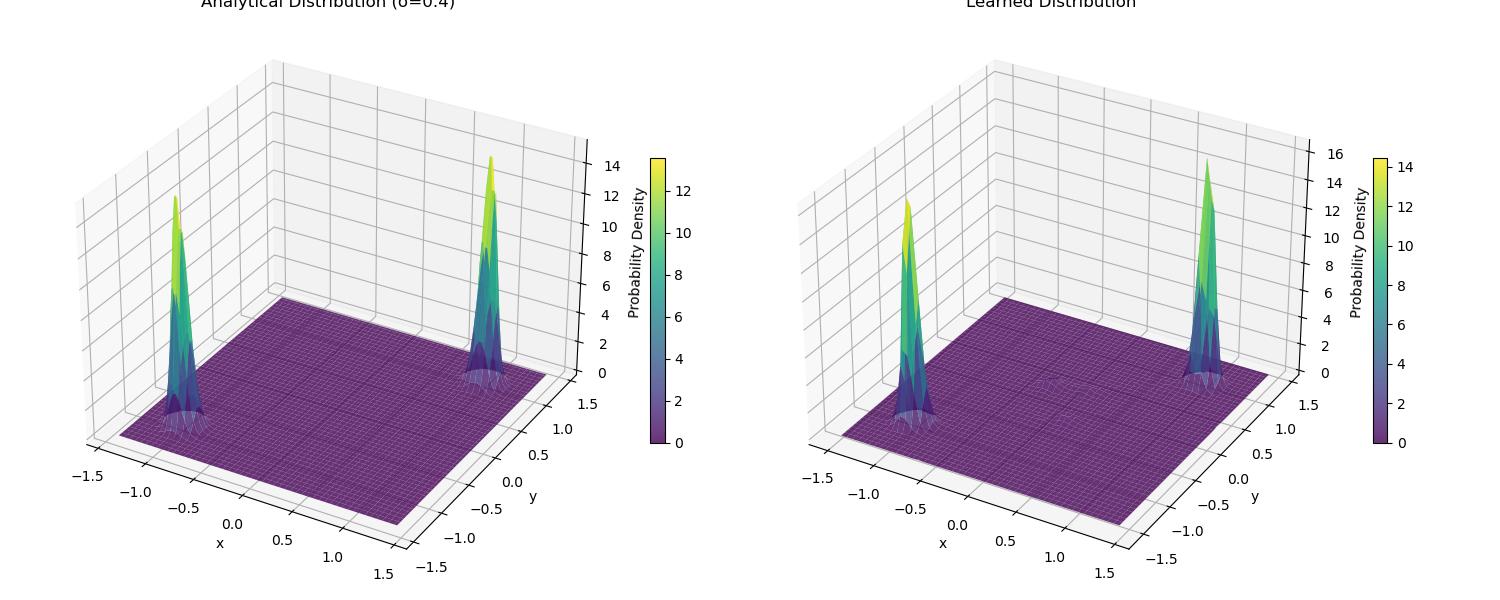}
    \caption{Comparison of analytical and learned probability distributions for the 2D double-peak problem. Left: Analytical Boltzmann distribution $\rho_{\text{true}}(\mathbf{x}) \propto \exp(-2V(\mathbf{x})/\sigma^2)$ with $\sigma=0.4$. Right: Learned distribution obtained from the pushforward network. The two distributions show excellent agreement in both peak locations and overall shape, validating the method's accuracy in capturing multimodal distributions.}
    \label{fig:2d_result_comparison}
\end{figure}

The training successfully converged to a solution that accurately captures the double-peak structure. As shown in Figure~\ref{fig:2d_result_training}, the scatter plot reveals clear bimodal clustering at the target positions $(-1, -1)$ and $(1, 1)$, with the density heatmap confirming probability concentration at these peaks. The loss convergence exhibits typical adversarial training dynamics, with rapid initial decrease followed by gradual decline and a transient instability around epochs 50,000--60,000 before stabilizing.

Figure~\ref{fig:2d_result_comparison} demonstrates excellent quantitative and qualitative agreement between the analytical and learned distributions, with sharp, well-defined peaks at correct locations and comparable peak heights and shapes.

\subsection{A Two-Dimensional Ring Distribution with Rotational Drift} \label{sec:ring}

We extend the framework to a ring-shaped distribution with divergence-free rotational drift, demonstrating the method's capability in handling non-gradient drift components and distributions with rotational symmetry.

\subsubsection{Problem Specification}

We consider the steady-state Fokker-Planck equation \eqref{eq:FPE} on $\mathbb{R}^2$ with a ring potential that has been smoothly extended to avoid singularities at the origin. The parameters are:

\begin{itemize}
    \item Dimension: $n=2$
    \item Potential: $V(x_1, x_2) = \frac{1}{4}(r^2 - r_0^2)^2$, where $r^2 = x_1^2 + x_2^2$
    \item Drift coefficient: $\mathbf{b}(\mathbf{x}) = \mathbf{b}_{\text{radial}} + \mathbf{b}_{\text{tangent}}$
    \item Diffusion coefficient: $a_{ij}(\mathbf{x}) = \sigma^2 \delta_{ij}/2$ (isotropic constant diffusion)
\end{itemize}

The potential $V(\mathbf{x})$ creates a smooth, ring-shaped energy landscape with a single minimum at radius $r = r_0$ and a local maximum at the origin. The drift has two distinct components:

\paragraph{Radial component (from the potential gradient):}
\begin{equation}
\mathbf{b}_{\text{radial}} = -\nabla V = -2(r^2 - r_0^2)\begin{pmatrix} x_1 \\ x_2 \end{pmatrix},
\end{equation}
which can be written explicitly as:
\begin{align}
b_1^{\text{rad}}(x_1, x_2) &= -2x_1(x_1^2 + x_2^2 - r_0^2), \\
b_2^{\text{rad}}(x_1, x_2) &= -2x_2(x_1^2 + x_2^2 - r_0^2).
\end{align}

\paragraph{Tangential component (rotation):}
\begin{equation}
\mathbf{b}_{\text{tangent}} = \omega \begin{pmatrix} -x_2 \\ x_1 \end{pmatrix},
\end{equation}
where $\omega$ is the angular velocity parameter controlling the rotation strength.

The tangential component is divergence-free:
\begin{equation}
\nabla \cdot \mathbf{b}_{\text{tangent}} = \frac{\partial(-\omega x_2)}{\partial x_1} + \frac{\partial(\omega x_1)}{\partial x_2} = 0,
\end{equation}
and therefore does not appear in the steady-state distribution, which follows the Boltzmann form:
\begin{equation}
\rho_{\text{true}}(\mathbf{x}) \propto \exp\left(-\frac{2V(\mathbf{x})}{\sigma^2}\right) = \exp\left(-\frac{(r^2 - r_0^2)^2}{2\sigma^2}\right).
\end{equation}

This distribution exhibits perfect rotational symmetry, with probability concentrated in a ring at radius $r_0$. The width of the ring is controlled by the noise parameter $\sigma$. The divergence-free drift $\mathbf{b}_{\text{tangent}}$ creates a steady probability current circulating around the ring but does not affect where probability accumulates.

For this test, we set $r_0 = 2$, $\omega = 2$, and $\sigma = 1$.

\subsubsection{Training Configuration}

The neural pushforward mapping $F_{\boldsymbol{\vartheta}}$ was implemented as a compact network with architecture $(4, 64, 2)$, using Tanh activation functions. The base distribution was a 4-dimensional standard Gaussian, $\mathbf{r} \sim \mathcal{N}(\mathbf{0}, \mathbf{I}_4)$.

The adversarial test functions were parameterized as $f^{(k)}(\mathbf{x}) = \sin(\mathbf{w}^{(k)} \cdot \mathbf{x} + b^{(k)})$, with $K=200$ test functions. Full adversarial training was implemented, alternating between generator updates (minimizing loss with respect to $\boldsymbol{\vartheta}$) and adversary updates (maximizing loss with respect to test function parameters), with 5 adversary updates per generator update every 10 epochs.

The model was trained for 20,000 epochs with batch size 1,000, using Adam optimizer with learning rate $10^{-2}$ for the generator and SGD with learning rate $10^{-2}$ for the test functions.

\subsubsection{Results and Performance}

\begin{figure}[htbp]
    \centering
    \includegraphics[width=\textwidth]{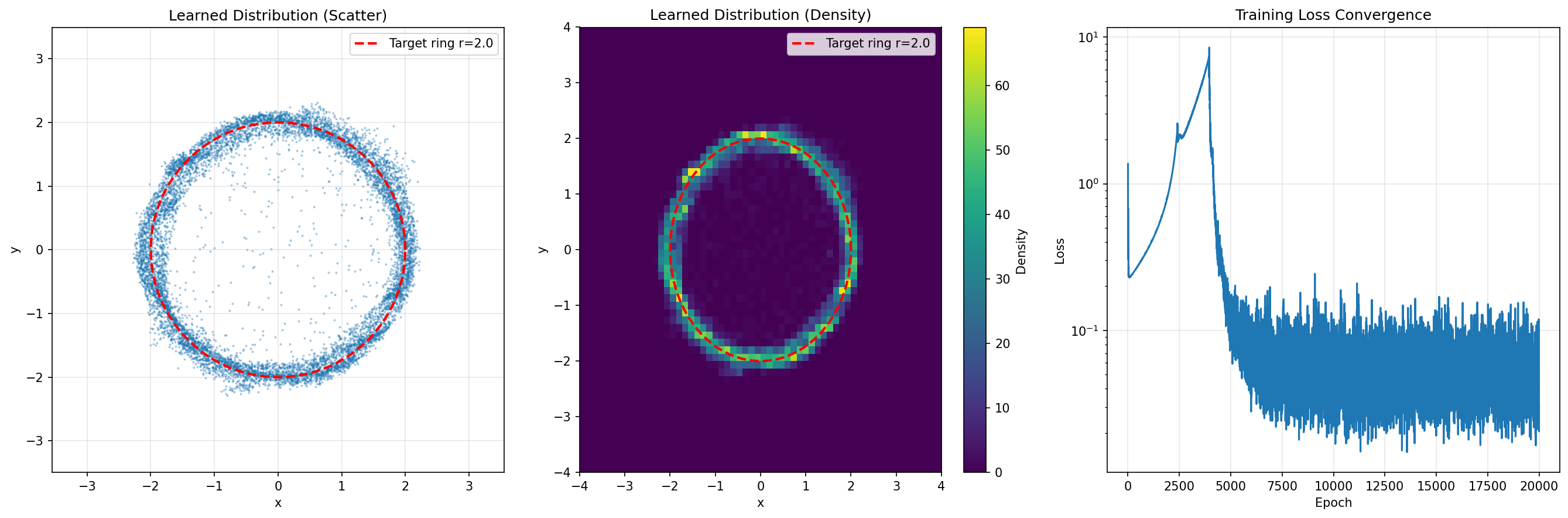}
    \caption{Training results for the two-dimensional ring potential with rotational drift. Left: Scatter plot of generated samples clearly showing concentration along the target ring at $r = r_0 = 2$ (red dashed circle). Middle: Density heatmap revealing the learned ring structure. Right: Training loss convergence showing initial adversarial dynamics followed by stable convergence.}
    \label{fig:ring_training}
\end{figure}

\begin{figure}[htbp]
    \centering
    \includegraphics[width=\textwidth]{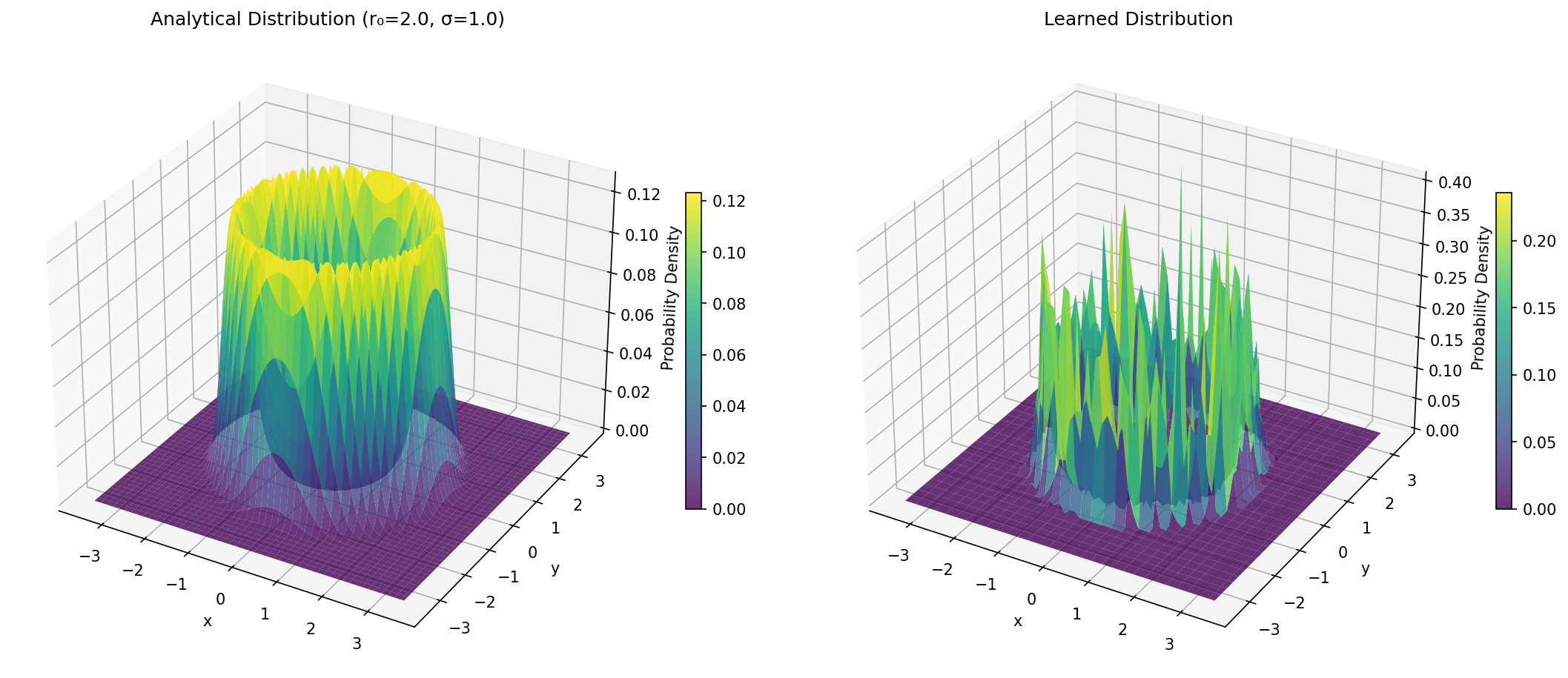}
    \caption{Comparison of analytical and learned probability distributions for the ring potential. Left: Analytical Boltzmann distribution showing smooth ring structure. Right: Learned distribution from the pushforward network. Both exhibit the characteristic ring shape at $r = r_0 = 2$.}
    \label{fig:ring_comparison}
\end{figure}

\begin{figure}[htbp]
    \centering
    \includegraphics[width=\textwidth]{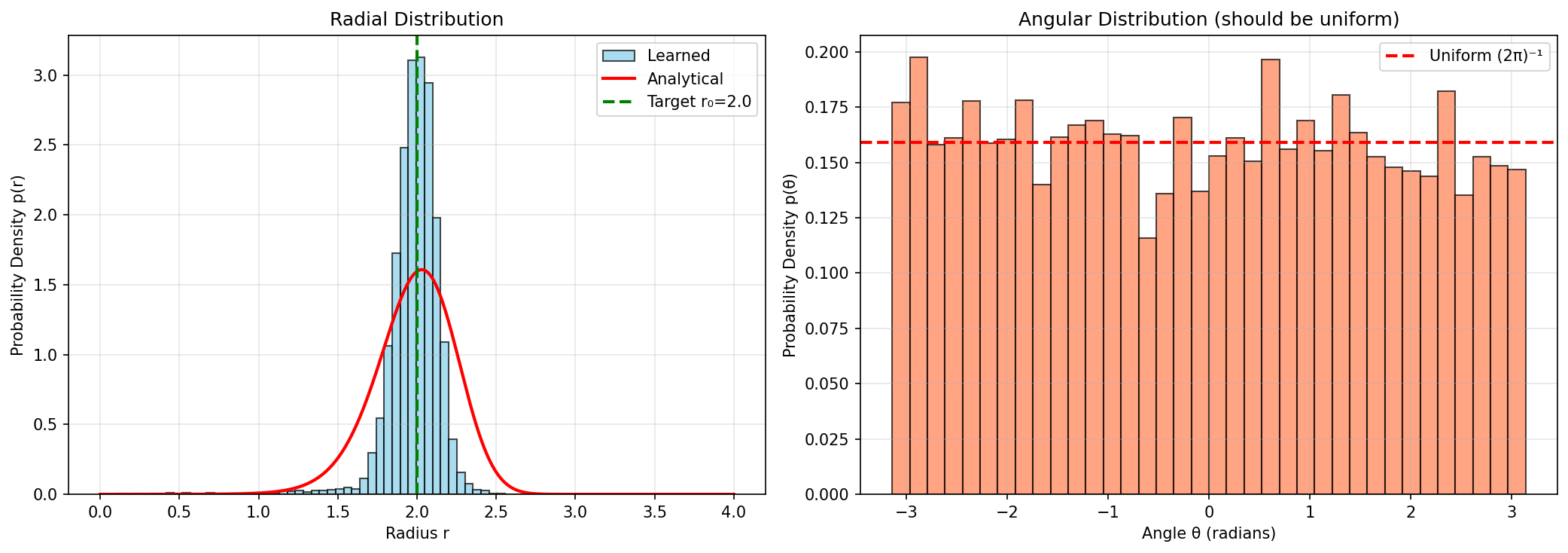}
    \caption{Decomposition of the learned ring distribution. Left: Radial probability density $p(r)$ showing agreement with the analytical prediction, with peak at $r = r_0 = 2$. Right: Angular probability density $p(\theta)$ demonstrating approximate uniformity, confirming the rotational symmetry of the steady-state distribution despite the presence of strong rotational drift ($\omega = 2$).}
    \label{fig:ring_radial_angular}
\end{figure}

\mbox{}
The training converged to a ring-shaped distribution, as evidenced by Figure~\ref{fig:ring_training}, with the loss decreasing from 1.36 to 0.032. Quantitative analysis reveals strong agreement with theoretical predictions: the mean radius $\langle r \rangle_{\text{learned}} = 1.979$ (target: $r_0 = 2.0$) and sample mean $\langle \mathbf{x} \rangle_{\text{learned}} = (-0.020, -0.004)$ (target: $\mathbf{0}$), with relative error approximately $1\%$ for the radius.

Figure~\ref{fig:ring_radial_angular} provides a detailed decomposition: the radial distribution shows a clear peak at $r = r_0$, while the angular distribution is approximately uniform, confirming that rotational symmetry is preserved despite the strong rotational drift. However, as shown in Figure~\ref{fig:ring_comparison}, the learned distribution exhibits a slightly higher peak density and narrower radial width compared to the analytical solution, suggesting that while the method correctly identifies the support and satisfies the weak form to reasonable accuracy (final loss $\approx 0.03$), achieving high-precision pointwise density estimates may require further refinement.

\subsection{Monte Carlo Integration using the Learned Sampler}

To demonstrate the practical utility of the learned pushforward network, we use it as a sampler for Monte Carlo integration and compare the results with high-accuracy grid-based quadrature. We take the neural poshforward mapping trained in \S\ref{sec:ring} consider computing expectations of the form
\begin{equation}
\mathbb{E}_{\rho}[f] = \int_{\mathbb{R}^2} f(x_1, x_2) \rho(x_1, x_2) \, dx_1 dx_2,
\end{equation}
for various test functions $f$, where $\rho$ is given as \eqref{ring-rho-true}.

The Monte Carlo estimate is obtained by generating $N$ samples $\{\mathbf{x}^{(i)}\}_{i=1}^N$ from the learned distribution via the pushforward network and computing
\begin{equation}
\mathbb{E}_{\rho}[f] \approx \frac{1}{N} \sum_{i=1}^N f(\mathbf{x}^{(i)}).
\end{equation}
For reference, we compute high-accuracy integrals using Simpson's rule on a $500 \times 500$ grid with the analytical Boltzmann distribution.

\begin{figure}[t!]
    \centering
    \includegraphics[width=\textwidth]{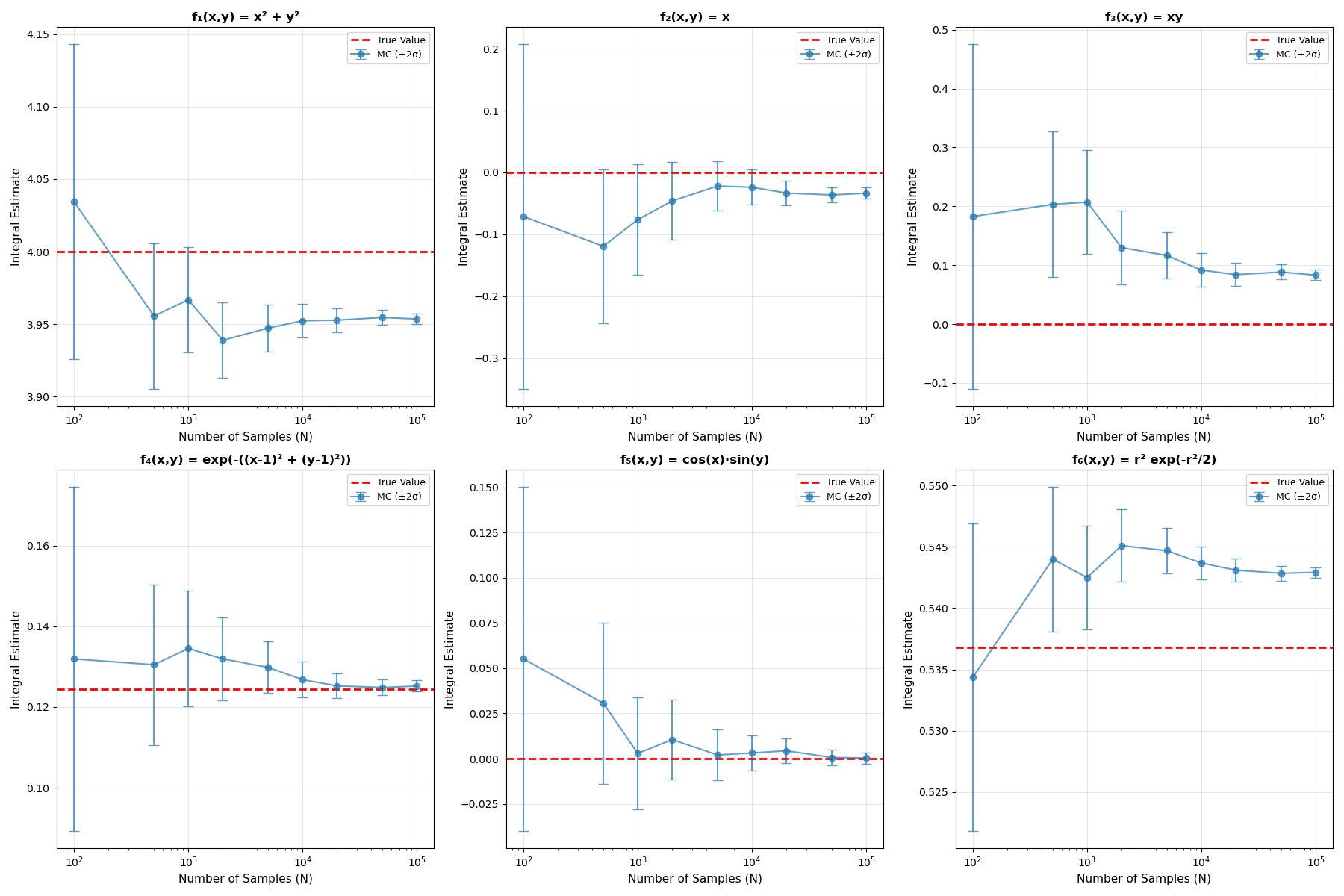}
    \caption{Monte Carlo integration convergence for six test functions using the learned ring distribution sampler. Blue curves with error bars show MC estimates ($\pm 2 \sigma$), and red dashed lines indicate grid-based reference value}
    \label{fig:mc_convergence}
\end{figure}

We test six different functions: $f_1(x,y) = x^2 + y^2$ (radial moment), $f_2(x,y) = x$ (first moment), $f_3(x,y) = xy$ (correlation), $f_4(x,y) = \exp(-(x-1)^2 - (y-1)^2)$ (localized Gaussian), $f_5(x,y) = \cos(x)\sin(y)$ (oscillatory), and $f_6(x,y) = r^2 \exp(-r^2)$ ($r = \sqrt{x^2+y^2}$. Figure~\ref{fig:mc_convergence} shows the convergence of Monte Carlo estimates as a function of sample size $N$.

The results demonstrate that the learned pushforward network serves as an effective importance sampler. For smooth functions like $f_4$ and $f_5$, the Monte Carlo estimates converge rapidly to the grid reference values, achieving errors below $10^{-3}$ with $10^4$ samples, with convergence rates closely following the expected $O(N^{-1/2})$ behavior.

\subsection{Time-Dependent One-Dimensional Ornstein-Uhlenbeck Process}

We now demonstrate the framework's capability for time-dependent Fokker-Planck equations by considering a one-dimensional Ornstein-Uhlenbeck process with Gaussian initial condition.

\subsubsection{Problem Specification}

We solve the time-dependent Fokker-Planck equation \eqref{eq:time_FPE} on $[0, T] \times \mathbb{R}$ with:
\begin{itemize}
    \item Dimension: $n = 1$
    \item Drift coefficient: $b(t, x) = -\theta(x - \mu)$, where $\theta = 1.0$ and $\mu = 0$
    \item Diffusion coefficient: $a(t, x) = \sigma^2/2$ (constant with $\sigma = 1.0$)
    \item Time interval: $T = 1.0$
    \item Initial condition: $\rho_0(x) = \mathcal{N}(\mu_0, \sigma_0^2)$, where $\mu_0 = 3.0$ and $\sigma_0 = 0.5$
\end{itemize}

The system starts from a Gaussian distribution positioned far from equilibrium and relaxes toward the steady state $\mathcal{N}(0, 0.5)$ over the time interval. The analytical solution for this initial value problem is:
\begin{align}
m(t) &= \mu_0 e^{-\theta t} + \mu(1 - e^{-\theta t}) = 3e^{-t}, \\
v(t) &= \sigma_0^2 e^{-2\theta t} + \frac{\sigma^2}{2\theta}(1 - e^{-2\theta t}) = 0.25e^{-2t} + 0.5(1 - e^{-2t}).
\end{align}

\subsubsection{Training Configuration}

For this time-dependent case, the initial condition term in the weak formulation is approximated via Monte Carlo:
\begin{equation}
\mathbb{E}_{x \sim \mathcal{N}(\mu_0, \sigma_0^2)}[f(0, x)] \approx \frac{1}{N} \sum_{i=1}^N f(0, x_i), \quad x_i \sim \mathcal{N}(\mu_0, \sigma_0^2),
\end{equation}
computed with $N = 500$ samples. The pushforward network architecture is:
\begin{equation}
F_\theta(t, r) = (\mu_0 + \sigma_0 \cdot z(r)) + t \cdot G_\theta(t, r),
\end{equation}
where $z(r)$ is a learned linear projection of the base noise $r \in \mathbb{R}^{d}$ to a standard normal variable, ensuring that $F_\theta(0, r) \sim \mathcal{N}(\mu_0, \sigma_0^2)$.

Training was performed with $K = 100$ test functions, batch size $M = 1000$, and base dimension $d = 2$ over 20,000 epochs.

\subsubsection{Results and Performance}

\begin{figure}[t!]
    \centering
    \includegraphics[width=\textwidth]{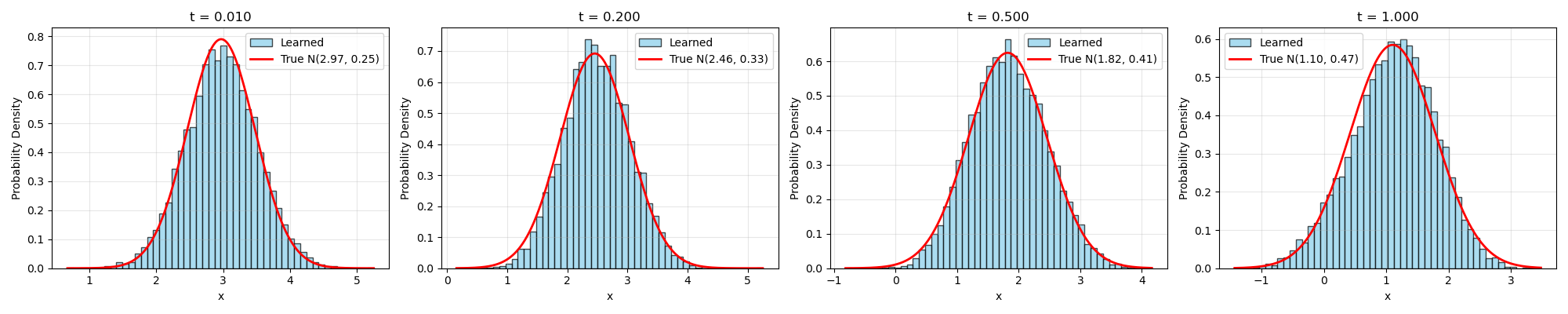}
    
    \vspace{0.5em}
    
    \begin{minipage}[t]{0.62\textwidth}
        \centering
        \includegraphics[width=\textwidth]{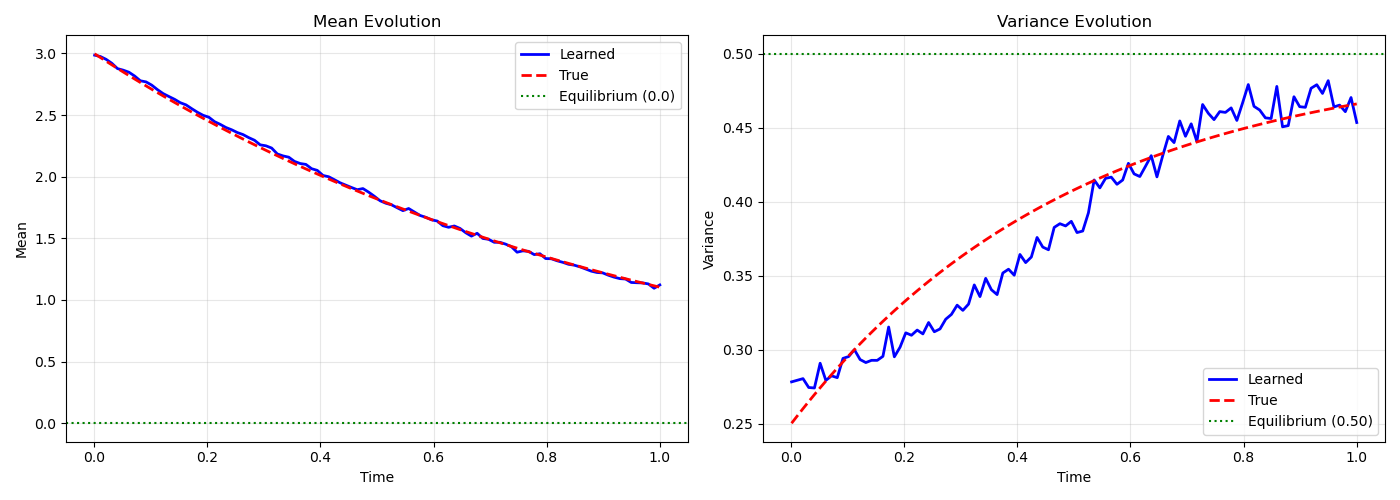}
    \end{minipage}
    \hfill
    \begin{minipage}[t]{0.37\textwidth}
        \centering
        \includegraphics[width=\textwidth]{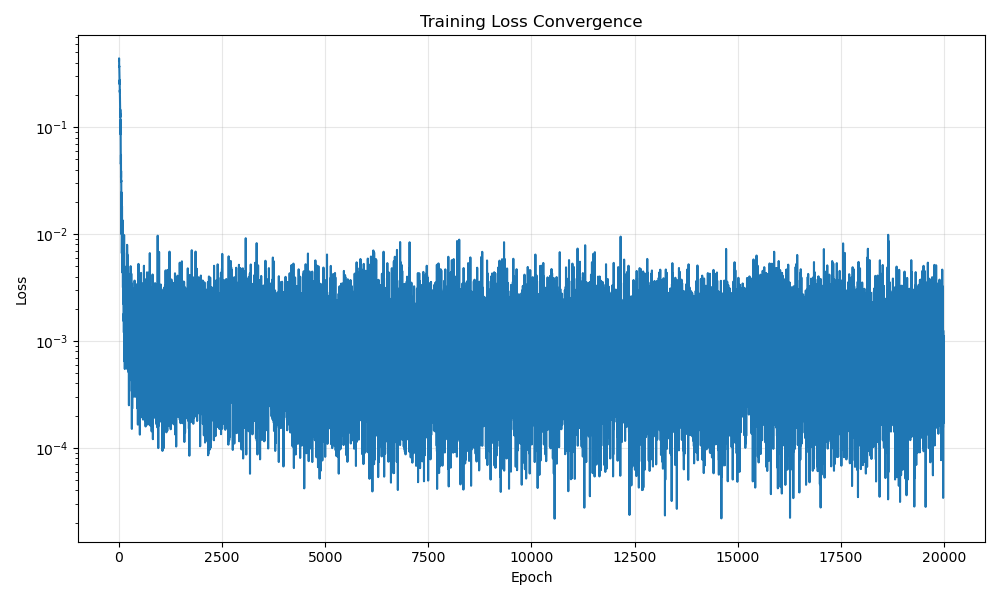}
    \end{minipage}
    
    \caption{Numerical results for the time-dependent Fokker-Planck equation with Gaussian initial condition $\mathcal{N}(3, 0.25)$ relaxing toward equilibrium $\mathcal{N}(0, 0.5)$. \textbf{Top:} Comparison of learned (histograms) and analytical (red curves) probability distributions at different times. \textbf{Bottom left\&{}middle:} Evolution of mean and variance over time \textbf{Bottom right:} Training loss convergence.}
    \label{fig:gaussian_ic_results}
\end{figure}

Figure~\ref{fig:gaussian_ic_results} demonstrates comprehensive validation of the learned solution. The top panel shows excellent agreement between learned and analytical distributions at times $t \in \{0.01, 0.2, 0.5, 1.0\}$, accurately capturing both the mean drift from $\mu_0 = 3$ toward equilibrium at $\mu = 0$ and the variance evolution from $\sigma_0^2 = 0.25$ toward the steady-state value of $0.5$. The bottom panels show the mean trajectory exhibiting perfect exponential decay with learned values tracking the analytical solution nearly exactly, while the variance shows characteristic relaxation with high accuracy despite some oscillations. The training loss converges from approximately $10^{-1}$ to below $10^{-4}$ within 1000 epochs, demonstrating rapid and stable convergence.

\subsection{High-Dimensional Time-Dependent Ornstein-Uhlenbeck Process}

We validate the framework's scalability by solving time-dependent Ornstein-Uhlenbeck processes in 10 and 100 dimensions. The problem specification is:
\begin{itemize}
\item Dimension: $n \in \{10, 100\}$
\item Drift: $\mathbf{b}(t, \mathbf{x}) = -\theta(\mathbf{x} - \boldsymbol{\mu})$ with $\theta = 1.0$, $\boldsymbol{\mu} = \mathbf{0}$
\item Diffusion: $a_{ij}(t, \mathbf{x}) = \frac{\sigma^2}{2}\delta_{ij}$ with $\sigma = 1.0$
\item Time interval: $T = 1.0$
\item Initial condition: $\rho_0(\mathbf{x}) = \mathcal{N}(\boldsymbol{\mu}_0, \sigma_0^2 \mathbf{I})$ with $\sigma_0 = 0.5$
\end{itemize}

For the 10D case, $\boldsymbol{\mu}_0 = (3.0, -2.5, 2.0, -3.5, 1.5, -1.0, 2.5, -2.0, 3.5, -1.5)^T$, placing each component far from equilibrium with different initial means. For the 100D case, initial means were drawn from $\mathcal{N}(0, 4)$ to ensure diverse starting positions across dimensions.

\subsubsection{Network Architecture and Training}

The pushforward network maps from $(t, x_0, \mathbf{r}) \in \mathbb{R}^{1+n+d}$ to $\mathbb{R}^n$ using the $\sqrt{t}$ scaling:
\begin{itemize}
\item \textbf{10D}: Base dimension $d=20$, hidden layers $[64, 64]$, trained for 20,000 epochs
\item \textbf{100D}: Base dimension $d=50$, hidden layers $[128, 128, 128]$, trained for 10,000 epochs
\end{itemize}

Training used batch sizes $M = 16{,}000$ (interior), $M_0 = M_T = 1{,}000$--$2{,}000$ (boundaries), with $K = 200$ test functions for 10D and $K = 5{,}000$ for 100D. The Adam optimizer was used for the generator ($\text{LR} = 5 \times 10^{-4}$ to $1 \times 10^{-3}$) and SGD for test functions ($\text{LR} = 1 \times 10^{-2}$), with adversarial updates every 10 epochs.

\subsubsection{Analytical Solution}

The analytical solution for each component $i$ evolves as:
\begin{align}
m_i(t) &= \mu_{0,i} e^{-\theta t}, \\
v(t) &= \sigma_0^2 e^{-2\theta t} + \frac{\sigma^2}{2\theta}(1 - e^{-2\theta t}),
\end{align}
with equilibrium distribution $\mathcal{N}(\mathbf{0}, \frac{\sigma^2}{2\theta}\mathbf{I}) = \mathcal{N}(\mathbf{0}, 0.5\mathbf{I})$.

\subsubsection{Results: 10-Dimensional Case}

\begin{figure}[htbp]
    \centering
    \includegraphics[width=0.7\textwidth]{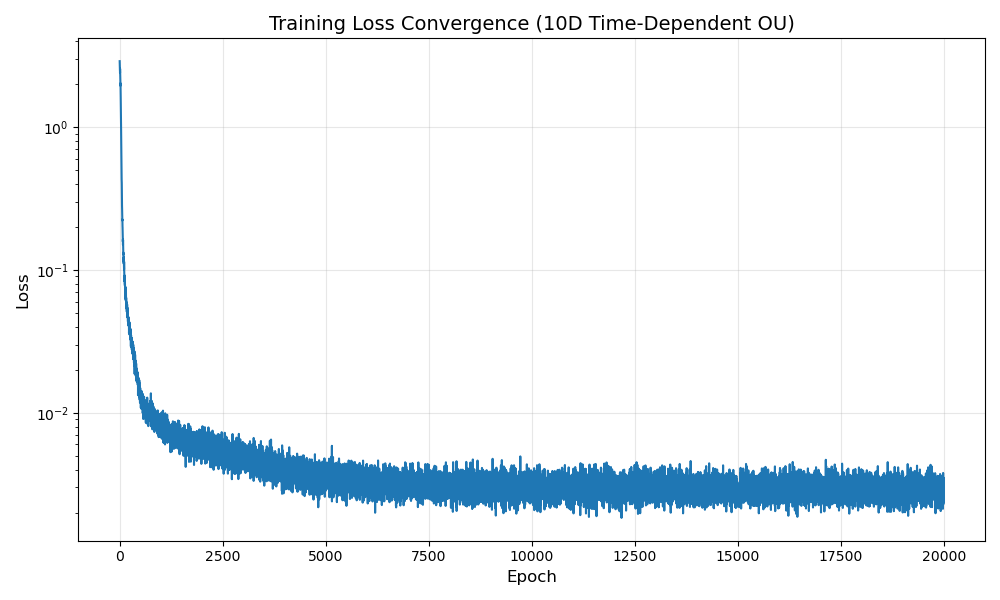}
    \caption{Training loss convergence for the 10D time-dependent OU process, demonstrating rapid decrease to below $10^{-3}$ within 1000 epochs.}
    \label{fig:10d_time_loss}
\end{figure}

\begin{figure}[htbp]
    \centering
    \includegraphics[width=\textwidth]{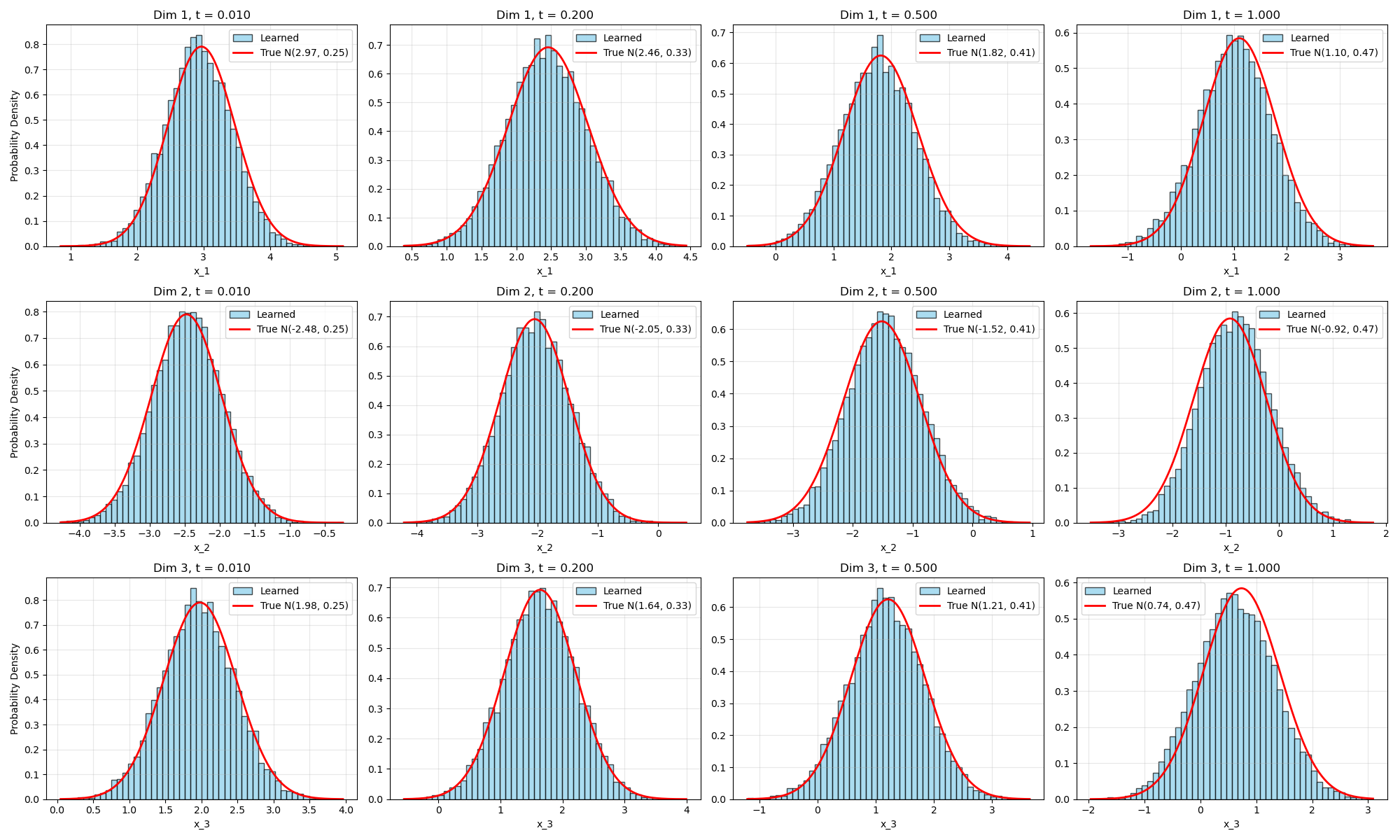}
    \caption{Probability distributions at different times for the 10D time-dependent OU process. Each row shows one spatial dimension at times $t \in \{0.01, 0.2, 0.5, 1.0\}$. Blue histograms represent learned distributions, while red curves show analytical Gaussian solutions.}
    \label{fig:10d_time_distributions}
\end{figure}

Figure~\ref{fig:10d_time_loss} shows rapid convergence from $10^0$ to below $10^{-3}$ within 1,000 epochs, then stabilizing around $10^{-3}$ to $10^{-4}$. The learned distributions (Figure~\ref{fig:10d_time_distributions}) accurately match analytical Gaussian solutions across all three displayed dimensions at times $t \in \{0.01, 0.2, 0.5, 1.0\}$, correctly capturing the relaxation dynamics from initial conditions toward equilibrium.

\begin{figure}[t!]
    \centering
    \includegraphics[width=\textwidth]{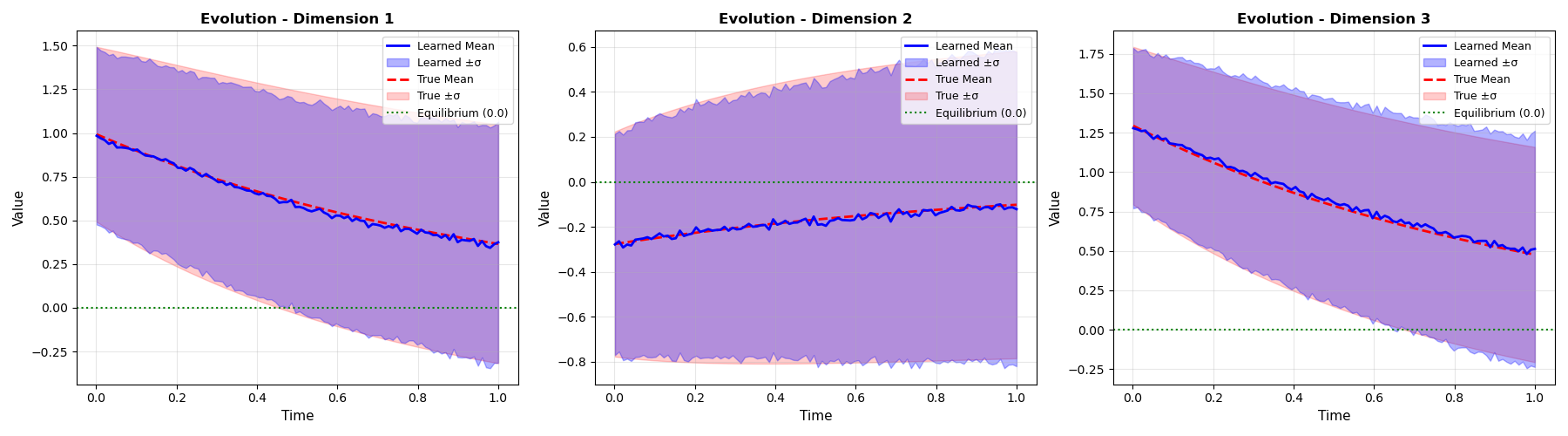}
    \caption{Temporal evolution of mean (top row) and variance (bottom row) for the first three dimensions of the 10D OU process. Blue solid lines show learned values, red dashed lines show analytical solutions, and green dotted lines indicate equilibrium values.}
    \label{fig:10d_time_evolution}
\end{figure}

The temporal evolution of first and second moments (Figure~\ref{fig:10d_time_evolution}) reveals excellent agreement for mean trajectories, which closely follow the analytical exponential decay $m_i(t) = \mu_{0,i} e^{-t}$ across all dimensions. Variance evolution shows characteristic relaxation from $\sigma_0^2 = 0.25$ toward equilibrium value $0.5$, matching the analytical solution $v(t) = 0.25 e^{-2t} + 0.5(1 - e^{-2t})$ despite some oscillations.

\subsubsection{Results: 100-Dimensional Case}

\begin{figure}[t!]
    \centering
    \includegraphics[width=0.7\textwidth]{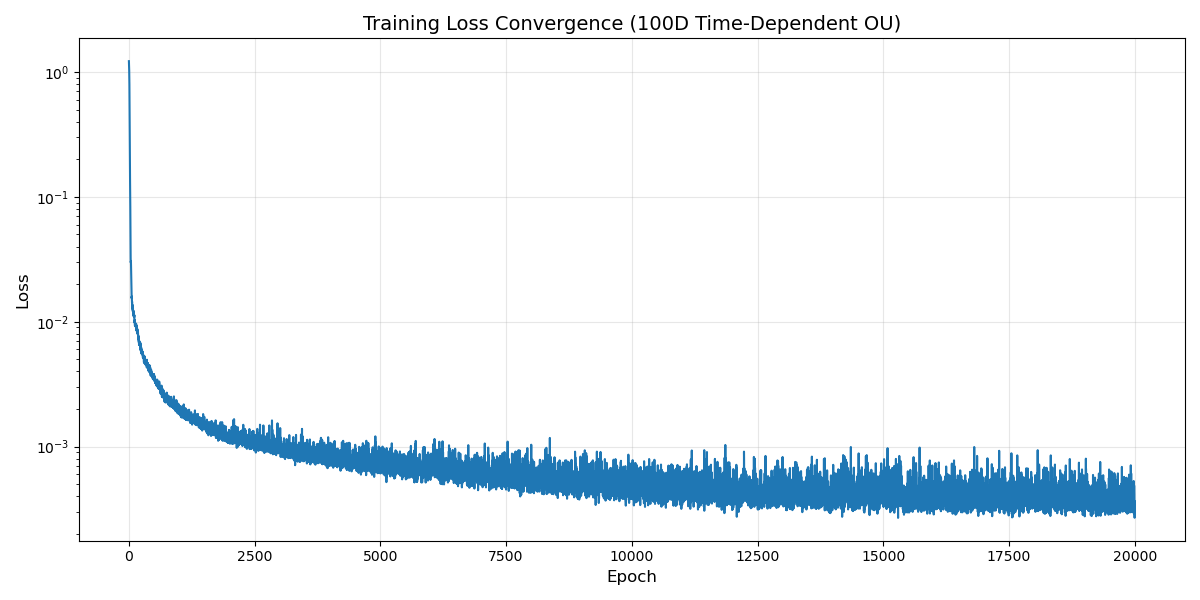}
    \caption{Training loss convergence for the 100D time-dependent OU process.}
    \label{fig:100d_time_loss}
\end{figure}

The 100D results demonstrate robust scalability. Training converged in 102 minutes on GPU, with loss decreasing from $10^0$ to below $10^{-3}$ within 2,500 epochs (Figure~\ref{fig:100d_time_loss}). Final statistics at $t=1.0$ show mean absolute error of 0.006 and variance error of 0.009, indicating high accuracy despite the dimensional complexity.

\begin{figure}[t!]
    \centering
    \includegraphics[width=\textwidth]{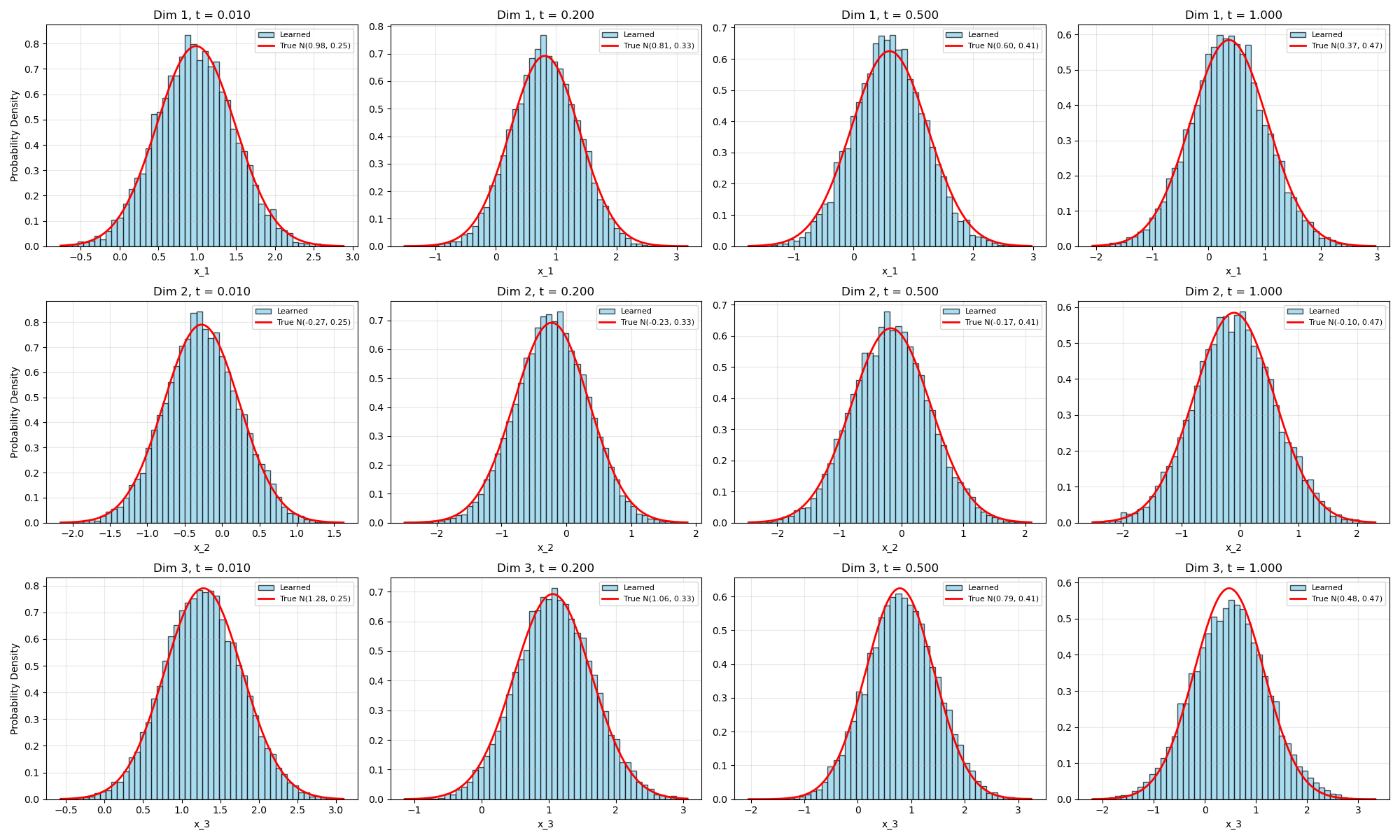}
    \caption{Probability distributions at different times for three representative dimensions of the 100D time-dependent OU process.}
    \label{fig:100d_time_distributions}
\end{figure}

\begin{figure}[t!]
    \centering
    \includegraphics[width=\textwidth]{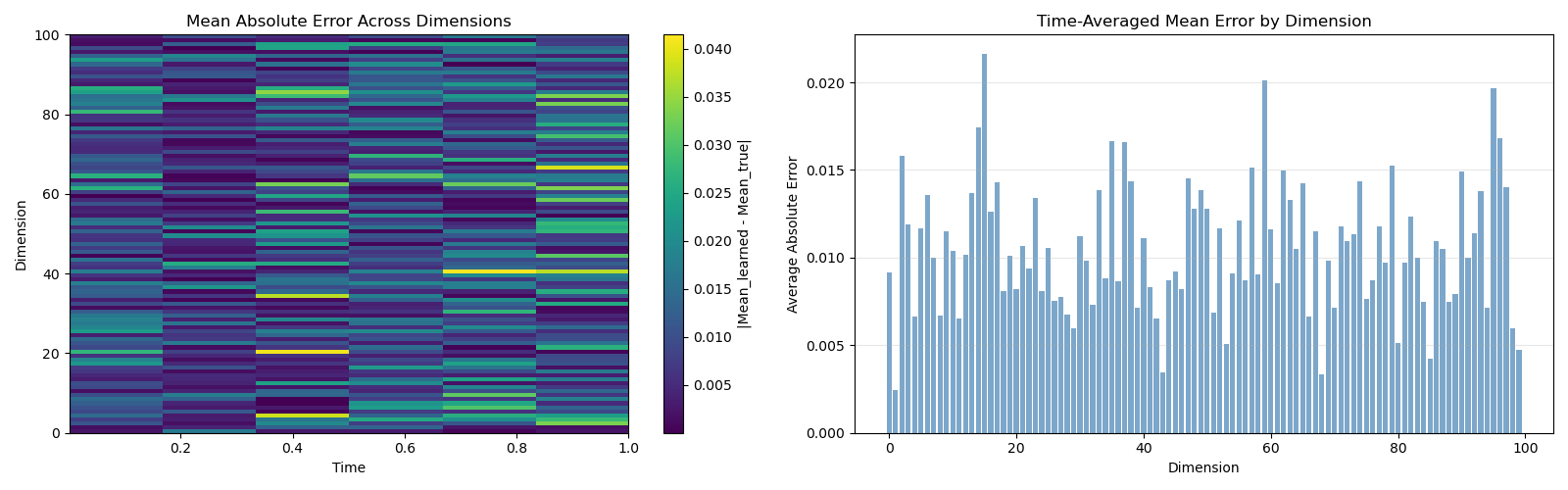}
    \caption{Error analysis for the 100D time-dependent OU process. Left: Heatmap of mean absolute errors across all 100 dimensions and validation time points. Right: Time-averaged mean error by dimension, showing consistent accuracy across the entire state space.}
    \label{fig:100d_error_heatmap}
\end{figure}

Figure~\ref{fig:100d_time_distributions} shows learned distributions for three representative dimensions at $t \in \{0.01, 0.2, 0.5, 1.0\}$, with excellent alignment to analytical Gaussians across all time points. The error heatmap (Figure~\ref{fig:100d_error_heatmap}) visualizes mean absolute errors across all 100 dimensions over time, revealing uniform accuracy with typical errors below 0.02 across most dimension-time combinations. Dimension-averaged errors show consistency across the entire state space, with no systematic bias toward particular dimensions.

\begin{figure}[t!]
    \centering
    \includegraphics[width=\textwidth]{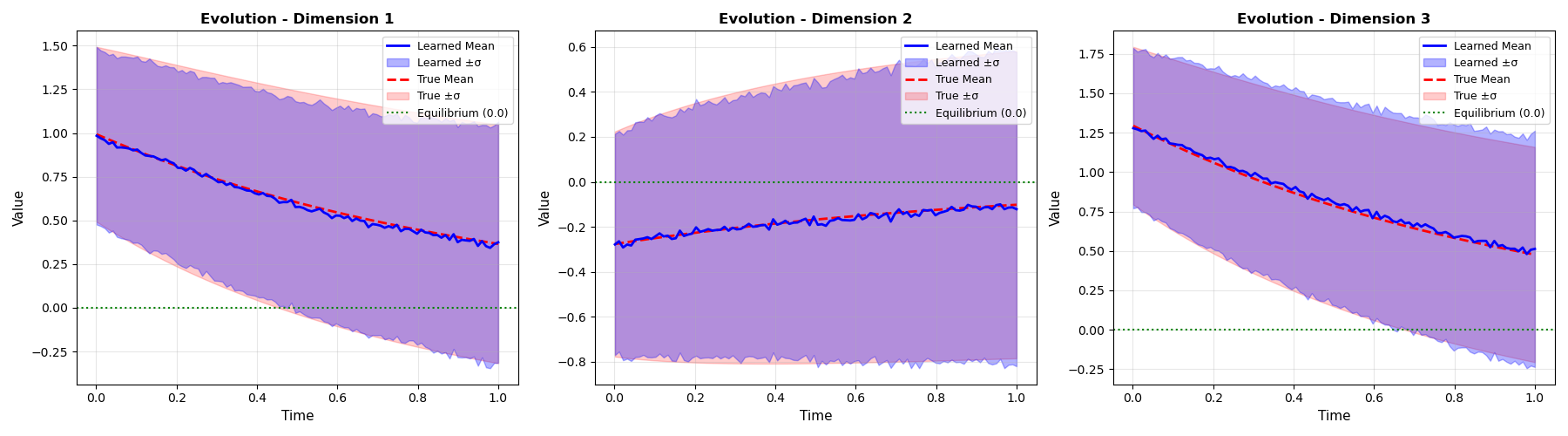}
    \caption{Temporal evolution of mean and variance for three representative dimensions of the 100D OU process, showing excellent agreement between learned (blue) and analytical (red) solutions.}
    \label{fig:100d_time_evolution}
\end{figure}

Temporal evolution plots (Figure~\ref{fig:100d_time_evolution}) for three sample dimensions confirm accurate mean tracking and variance relaxation, matching analytical predictions. The mean evolution captures the exponential decay from diverse initial conditions toward equilibrium, while variance evolution correctly interpolates between initial spread ($\sigma_0^2 = 0.25$) and equilibrium variance ($0.5$).

These high-dimensional results validate the framework's capability to handle complex time-dependent problems with many degrees of freedom. The method maintains accuracy while scaling from 10 to 100 dimensions, requiring only modest increases in network capacity and test function count. The $\sqrt{t}$ parameterization effectively handles the diffusive time scaling, and the three-term weak formulation correctly enforces both initial conditions and interior dynamics.

\section{Conclusion}

We have introduced a weak adversarial framework for solving Fokker-Planck equations by learning neural pushforward samplers. The method transforms the PDE solution problem into learning a pushforward map from a simple base distribution to the target solution distribution, with training guided by a weak formulation that permits arbitrary network architectures and employs computationally efficient plane-wave test functions.

Our numerical experiments validate the framework across diverse problem settings. In steady-state problems, the method successfully learned unimodal Gaussian distributions (1D Ornstein-Uhlenbeck), complex bimodal distributions arising from double-well potentials (1D and 2D), and ring distributions with non-gradient rotational drift components. The learned samplers demonstrated practical utility for Monte Carlo integration, achieving convergence rates consistent with theoretical expectations. For time-dependent problems, the framework accurately captured the relaxation dynamics of Ornstein-Uhlenbeck processes in both one and three dimensions, with learned solutions showing excellent agreement with analytical mean and variance evolution.

The results demonstrate that plane-wave test functions, despite their simple form, provide sufficient expressiveness to capture complex distributional features through adversarial optimization. The framework naturally handles multimodal distributions and probability currents induced by divergence-free drift, while maintaining probability conservation without explicit normalization constraints. The extension from steady-state to time-dependent problems and from one to three dimensions proceeded straightforwardly, requiring only appropriate increases in test function count and network capacity.

The primary limitation observed is that while the method accurately captures distribution support and moment statistics, achieving high-precision pointwise density estimates may require refinement. This is expected given the weak formulation's focus on integrated conditions rather than pointwise values. Future directions include extending to higher dimensions where traditional grid-based methods become prohibitive, adaptive selection of test functions, and hybrid approaches combining weak form satisfaction with targeted density matching in regions of interest.



\begin{thebibliography}{1}
\bibitem{PINN}
M.~Raissi, P.~Perdikaris, and G.~E. Karniadakis.
\newblock Physics-informed neural networks: A deep learning framework for solving forward and inverse problems involving nonlinear partial differential equations.
\newblock {\em Journal of Computational Physics}, 378:686--707, 2019.
\newblock doi: 10.1016/j.jcp.2018.10.045.

\bibitem{socmartnet}
Wei Cai, Shuixin Fang, and Tao Zhou.
\newblock {SOC-MartNet: A Martingale Neural Network for the Hamilton-Jacobi-Bellman Equation without Explicit inf H in Stochastic Optimal Controls}, 2025.

\bibitem{dfsocmartnet}
Wei Cai, Shuixin Fang, Tao Zhou, and Wenzhong Zhang.
\newblock {A Derivative-Free Martingale Neural Network SOC-Martnet For The Hamilton-Jacobi-Bellman Equations In Stochastic Optimal Controls}, 08 2024.

\bibitem{Chen2025}
Chen Chen, Yunan Yang, Yang Xiang, and Others.
\newblock Automatic differentiation is essential in training neural networks for solving differential equations.
\newblock {\em Journal of Scientific Computing}, 104(54), 2025.

\bibitem{weran}
Zhizhong Kong, Rui Sheng, Jerry~Zhijian Yang, Juntao You, and Cheng Yuan.
\newblock {A Weaker Adversarial Neural Networks Method for Linear Second Order Elliptic PDEs}, 2025.
\newblock Available at SSRN.

\bibitem{wan}
Yaohua Zang, Gang Bao, Xiaojing Ye, and Haomin Zhou.
\newblock Weak adversarial networks for high-dimensional partial differential equations.
\newblock {\em Journal of Computational Physics}, 411:109409, 2020.

\bibitem{drdm}
W.~Cai, S.~Fang, and T.~Zhou.
\newblock Deep random difference method for high dimensional quasilinear parabolic partial differential equations.
\newblock {\em Communications in Mathematics and Statistics}, 11(3):599--626, 2023.
\newblock doi: 10.1007/s40304-021-00263-6.

\bibitem{deepritz}
Weinan E and Bing Yu.
\newblock The Deep Ritz Method: A Deep Learning-Based Numerical Algorithm for Solving Variational Problems.
\newblock {\em Communications in Mathematics and Statistics}, 6(1):1--12, March 2018.

\bibitem{NPFP}
\newblock Shu Liu, Wuchen Li, Hongyuan Zha, and Haomin Zhou
\newblock Neural Parametric Fokker--Planck Equation.
\newblock {\em SIAM Journal on Numerical Analysias} 2022 60:3, 1385-1449

\bibitem{WHF}
\newblock Shui-Nee Chow, Wuchen Li, Haomin Zhou,
\newblock Wasserstein Hamiltonian flows,
\newblock {\em Journal of Differential Equations},
Volume 268, Issue 3,
2020,
Pages 1205-1219,
ISSN 0022-0396,
https://doi.org/10.1016/j.jde.2019.08.046.

\bibitem{XWanFP}
\newblock Kejun Tang, Xiaoliang Wan, Qifeng Liao,
\newblock Adaptive deep density approximation for Fokker-Planck equations,
\newblock {\em Journal of Computational Physics},
Volume 457,
2022,
111080,
ISSN 0021-9991,
https://doi.org/10.1016/j.jcp.2022.111080.

\bibitem{TZhouFP}
\newblock Xiaodong Feng, Li Zeng, Tao Zhou, 
\newblock Solving Time Dependent Fokker-Planck Equations via Temporal Normalizing Flow, 
\newblock {\em Communications in Computational Physics}, Vol. 32 (2022), Iss. 2 : pp. 401–423

\bibitem{realNVP}
Laurent Dinh, Jascha Sohl-Dickstein, and Samy Bengio.
\newblock Density estimation using Real NVP.
\newblock \emph{International Conference on Learning Representations (ICLR)}, 2017.
\newblock arXiv:1605.08803.

\bibitem{GAN}
Ian Goodfellow, Jean Pouget-Abadie, Mehdi Mirza, Bing Xu, David Warde-Farley, Sherjil Ozair, Aaron Courville, and Yoshua Bengio.
\newblock Generative adversarial nets.
\newblock In \emph{Advances in Neural Information Processing Systems}, pages 2672--2680, 2014.

\end{thebibliography}
\end{document}